\documentclass[leqno,11pt]{amsart}

\usepackage{amssymb, amsmath}
\def\refer#1{~\ref{#1}}
\def\refeq#1{~(\ref{#1})}
\def\ccite#1{~\cite{#1}}

\def\suite#1#2#3{(#1_{#2})_{#2\in {#3}}}

\def\longformule#1#2{
\displaylines{ \qquad{#1} \hfill\cr \hfill {#2} \qquad\cr } }
\def\inte#1{
\displaystyle\mathop{#1\kern0pt}^\circ }

\def\supetage#1#2{
\sup_{\scriptstyle {#1}\atop\scriptstyle {#2}} }

\let\al=\alpha

\let\e=\varepsilon

\let\lam=\lambda

\let\D=\Delta

\let\wt=\widetilde
\let\wh=\widehat

\def\cB{{\mathcal B}}
\def\cC{{\mathcal C}}

\def\cF{{\mathcal F}}

\def\cS{{\mathcal S}}
\def\cT{{\mathcal T}}

\def\virgp{\raise 2pt\hbox{,}}
\def\cdotpv{\raise 2pt\hbox{;}}

\def\eqdefa{\buildrel\hbox{{\rm \footnotesize def}}\over =}
\def\Id{\mathop{\rm Id}\nolimits}

\def\C{\mathop{\mathbb C\kern 0pt}\nolimits}
\def\DD{\mathop{\mathbb D\kern 0pt}\nolimits}
\def\EE{\mathop{{\mathbb E \kern 0pt}_\e}\nolimits}
\def\K{\mathop{\mathbb K\kern 0pt}\nolimits}
\def\N{\mathop{\mathbb N\kern 0pt}\nolimits}
\def\Q{\mathop{\mathbb Q\kern 0pt}\nolimits}
\def\R{\mathop{\mathbb R\kern 0pt}\nolimits}
\def\SS{\mathop{\mathbb S\kern 0pt}\nolimits}
\def\ZZ{\mathop{\mathbb Z\kern 0pt}\nolimits}
\def\TT{\mathop{\mathbb T\kern 0pt}\nolimits}
\def\P{\mathop{\mathbb P\kern 0pt}\nolimits}

\newcommand{\ds}{\displaystyle}

\def\dive{\mathop{\rm div}\nolimits}

\newcommand{\andf}{\quad\hbox{and}\quad}
\newcommand{\with}{\quad\hbox{with}\quad}
\newcommand{\beq}{\begin{equation}}
\newcommand{\eeq}{\end{equation}}
\newcommand{\ben}{\begin{eqnarray}}
\newcommand{\een}{\end{eqnarray}}
\newcommand{\beno}{\begin{eqnarray*}}
\newcommand{\eeno}{\end{eqnarray*}}

\def\vbar{\overline v}

\newtheorem{defi}{Definition}[section]
\newtheorem{theo}{Theorem}
\newtheorem{lemma}{Lemma}[section]

\newtheorem{corol}{Corollary}[section]
\newtheorem{prop}{Proposition}[section]

\begin{document}
%%%%%%%%%%%%%%%%%%%%%%%%%%%%%%%%%%%%%%%%%%%%%%%%%%%%%%%%%%%%%%%%%%%%%%
%%  Text                                                            %%
%%%%%%%%%%%%%%%%%%%%%%%%%%%%%%%%%%%%%%%%%%%%%%%%%%%%%%%%%%%%%%%%%%%%%%

\title[Large solutions to the Navier-Stokes equations]{Global regularity for  some classes of large solutions    to
  the   Navier-Stokes equations}

\author[J.-Y. Chemin]{Jean-Yves  Chemin}
\address[J.-Y. Chemin]%
{ Laboratoire J.-L. Lions UMR 7598\\ Universit{\'e} Pierre et Marie Curie\\
175, rue du Chevaleret\\ 75013 Paris\\FRANCE }
\email{chemin@ann.jussieu.fr }
\author[I. Gallagher]{Isabelle Gallagher}
\address[I. Gallagher]%
{ Institut de Math{\'e}matiques de Jussieu UMR 7586\\ Universit{\'e} Paris Diderot\\
175, rue du Chevaleret\\ 75013 Paris\\FRANCE }
\email{Isabelle.Gallagher@math.jussieu.fr}
\author[M. Paicu]{Marius Paicu}
\address[M. Paicu]{Universit\'e Paris Sud\\Laboratoire de Math\'ematiques\\
B\^atiment 425\\
91405 Orsay Cedex\\ France}
\email{Marius.Paicu@math.u-psud.fr}

%\date{}

\begin{abstract} In \cite{cgens}-\cite{cg3} classes of initial data to the three dimensional, incompressible Navier-Stokes equations were presented, generating a global smooth solution although the norm of the initial data may be chosen arbitrarily large. 
The main feature of the initial data considered in~\cite{cg3} is that it varies slowly in one direction, though in some sense it is ``well prepared'' (its norm is large but does not depend on the slow parameter).Ê
 The aim of this article is to generalize the setting of~\cite{cg3} to an ``ill prepared'' situation (the norm blows up as the small parameter goes to zero).
 As in~\cite{cgens}-\cite{cg3}, the proof uses the special structure of the nonlinear term of the equation.
\end{abstract}


\maketitle

%%%%%%%%%%%%%%% %%%%%%% %%% INTRODUCTION %%%%%%%%%%%%%%%%%%%%%

\setcounter{equation}{0}
\section{Introduction}\label{intro}
\setcounter{equation}{0} 
 
We study in this paper the Navier-Stokes equation with initial data which are slowly varying in the vertical variable. More precisely we consider the system
$$
\begin{cases}
\partial_t u+u\cdot\nabla u- \Delta u=-\nabla p \quad\text{in}\quad \R^+\times \Omega \\
\dive u=0\\
u|_{t=0}=u_{0,\e},
\end{cases}
$$
where $\Omega=\TT^2\times\R$ (the choice of this particular domain will be explained later on) and $u_{0,\e}$  is a divergence free vector field, whose dependence on the vertical variable~$x_3$ will be chosen to be ``slow'', meaning that it depends on~$\e x_3$ where~$\e$ is a small parameter.    Our goal is to prove a global existence in time result for the solution generated by this type of   initial data, with no smallness assumption on its norm.

  \subsection{Recollection of some known results on the Navier-Stokes equations}

The mathematical study of the Navier-Stokes equations has a long history, which we shall describe briefly in this paragraph. We shall first recall results concerning the main global wellposedness results, and some blow-up criteria. Then we shall concentrate on the case when the special algebraic structure of the system is used, in order to improve those previous results. 

To simplify we shall place ourselves in the whole euclidian space~$\R^d$ or in the torus~$\TT^d$ (or in variants of those spaces, such as~$\TT^2 \times \R$ in three space dimensions); of course results exist in the case when the equations are posed in domains of the euclidian space, with Dirichlet boundary conditions, but we choose to simplify the presentation by not   mentioning explicitly those studies (although some of  the theorems recalled below also hold in the case of domains up to obvious modifications of the statements and sometimes much more difficult proofs).

\subsubsection{Global wellposedness and blow-up results}
 The first important result on the  Navier-Stokes system
was obtained by J. Leray  in his seminal paper\ccite{leray} in 1933. He proved
that any finite energy initial data (meaning square-integrable data) generates a (possibly non unique) global in time weak
solution which satisfies an energy estimate; he moreover proved in~\cite{leray2D} the uniqueness of the solution in two space dimensions. Those results use the structure of the nonlinear terms, in order to obtain the energy inequality. He also proved the uniqueness of   weak solutions in three space dimensions, under the additional condition that one of the weak solutions has more regularity properties (say belongs to~$L^2(\R^+;L^\infty)$:  this would now be  qualified as a ``weak-strong uniqueness result''). The question of the global wellposedness of the Navier-Stokes equations was then raised, and has been open ever since. We shall now  present a few of the historical landmarks in that study.

The Fujita-Kato theorem \cite{fujitakato}   gives a partial answer
to the construction of a global unique solution. Indeed, that theorem provides a unique,  local in time solution in the homogeneous Sobolev space   $\dot H^{\frac d2-1}$ in $d$ space dimensions, and that solution is proved to be global if the initial data is small 
in~$\dot H^{\frac d2-1}$ (compared to the viscosity, which is chosen equal to one here to simplify). The result was improved to the Lebesgue space $L^d $  by F. Weissler in\ccite{weisslerNS} (see also \ccite{gigamiyakawa} and\ccite{kato}).  The  method   consists in applying a Banach fixed point theorem to the integral formulation of the equation, and was generalized by M. Cannone, Y. Meyer and  F. Planchon in~\cite{cannonemeyerplanchon} to Besov spaces of negative index of regularity. More precisely they proved that  if the initial data is small in the Besov space
$\dot B^{-1+\frac dp}_{p,\infty}$ (for~$p<\infty$), then there is a unique,  global in time solution. Let us emphasize that this result
allows to construct global solutions for strongly oscillating initial data
which may have a large norm in $\dot H^{\frac d2-1} $ or in $L^d$. A typical example in three space dimensions is
$$
u_0^\e(x)Ê\eqdefa \e^{-\alpha}\sin\Bigl(\frac{x_3}{\e}\Bigr)\bigl(-\partial_2\varphi(x),\partial_1\varphi(x),0\bigr),
$$
where $0<\alpha<1$ and $\varphi\in\mathcal S(\R^3;\R)$. This can be checked by using the definition of Besov norms:
$$
\forall s >0, \:\forall (p,q) \in [1,\infty], \quad \|f\|_{\dot B^{-s}_{p,q}} \eqdefa \left\| t^\frac s2 \|e^{t\Delta} f\|_{L^p} \right\|_{L^q(\R^+;\frac{dt}t)}.
$$
 More recently  in \cite{kochtataru}, H. Koch and D. Tataru obtained a unique global in time
solution for data small enough in a more general space, consisting of vector fields whose components are derivatives of $\text{BMO}$
functions. The norm in that space is given by
\beq \label{defBMO-1}
\|u_0\|^2_{BMO^{-1}}\eqdefa \sup_{t>0} t \| e^{t\Delta}u_0\|^2_{L^\infty}
+\supetage {x\in \R^d}{R>0} \frac 1 {R^d} \int_{P(x,R)} |(e^{t\Delta} u_0)(t,y)|^2dy,
\eeq
where $P(x,R)$ stands for the parabolic set~$[0,R^{2}] \times B(x,R)$ while~$B(x,R)$ is the ball centered at~$x$, of radius~$R.$
  
  \smallskip

One should notice that spaces where global, unique solutions are constructed  for small initial data, are necessarily scaling-invariant spaces: thus all those spaces are invariant under the invariant transformation for Navier-Stokes equation  $u_\lambda(t,x)=\lambda u(\lambda^2 t,\lambda x)$. Moreover it can be proved  (as observed for instance in\ccite{cg2}) that if~$B$ is  a Banach space
continuously included in the space~$\cS'$ of tempered distributions  such  
that
 $$
\mbox{for  any} \: (\lam,  a)\in \R^+_{\star}\times \R^d , Ê\quad \|f(\lam(\cdot-  a))\|_{B} 
 = \lam^{-1}\|f\|_{B}, $$
then $
 \|\cdot\|_B\leq C \ds \sup_{t>0}  t^{\frac 1 2} \| e^{t\Delta}u_0\|_{L^\infty}.
 $
 One recognizes on the right-hand side of the inequality the~$\dot B^{-1 }_{\infty,\infty}$ norm, which is slightly smaller than the~$BMO^{-1}$ norm recalled above in~(\ref{defBMO-1}): indeed  the~$BMO^{-1}$ norm takes into account not only the~$\dot  B^{-1 }_{\infty,\infty}$ information,  but also the fact that  first Picard iterate of the Navier-Stokes equations should be locally square integrable in space and time. It thus seems that the Koch-Tataru theorem is optimal for the wellposedness of the Navier-Stokes equations. This observation also shows that if one wants to go beyond a smallness assumption on the initial data to prove the global existence of unique solutions, one should check that the~$\dot B^{-1 }_{\infty,\infty}$ norm of the initial data may be chosen large.

 To conclude this paragraph, let us remark that the fixed-point methods used to prove local in time wellposedness for arbitrarily large data (such results are available in Banach spaces in which the Schwartz class is dense, typically~$\dot B^{-1+\frac dp}_{p,q}$ for finite~$p$ and~$q$) naturally provides blow-up criteria. For instance, one can prove that if the life span of the solution is finite, then  the~$L^q([0,T];\dot B^{-1+\frac dp + \frac2q}_{p,q})$ norm blows up as~$T$ approaches the blow up time. A natural question is to ask if the~$\dot B^{-1+\frac dp}_{p,q}$ norm itself blows up. Progress has been made very recently on this question, and uses the specific structure of the equation, which was not the case for the results presented in this paragraph. We therefore postpone the exposition of those results to the next paragraph.
 
 We shall not describe more results on the Cauchy problem for the Navier-Stokes equations, but refer the interested reader for instance to the monographs\ccite{lemarie} and\ccite{meyer} for more details.
 
%
%Concerning the methods to obtain such results, We recall that proving the existence of a unique, global in time solution to the Navier-Stokes equation is
% rather standard  (it is a consequence of a Banach fixed point theorem) as long as the initial data is chosen small enough in some scale invariant spaces embedded in~$\dot C^{-1}$, where we recall that
%$$
%\|f\|_{\dot C^{-1}} \eqdefa \sup_{t > 0} t^\frac12 \|e^{t\Delta} f\|_{L^\infty}.
%$$
%We refer for instance to~\cite{cannonemeyerplanchon},\cite{fujitakato},\cite{gigamiyakawa},\cite{kochtataru},\cite{weisslerNS} for a proof in various scale invariant function spaces.  This theorems are generals results of global existence for small initial data and does not take into account the any particular  algebrical properties of the nonlinear terms in the Navier-Stokes equation

  \subsubsection{Results using the specific algebraic structure of the equation}
  If one wishes to improve the theory on the Cauchy problem for the Navier-Stokes equations, it seems crucial to use the specific structure of the nonlinear term in the equations, as well as the divergence-free assumption. Indeed it was proved by S. Montgomery-Smith in\ccite{ms} (in a one-dimensional setting, which was later generalized to a 2D and 3D situation by two of the authors  in~\cite{gp}) that some models exist for which  finite time blow up can be proved for some classes of large data, despite the fact that the same  small-data global wellposedness results hold as for the Navier-Stokes system. 
  Furthermore, the generalization to the 3D case in\ccite{gp} shows that some large initial data which generate  a global solution for the Navier-Stokes equations (namely the data of~\cite{cgens} which will be presented below) actually generate a blowing up solution for the toy model. 
  
In this paragraph, we shall present a number of wellposedness theorems (or blow up criteria) which have been obtained in the past and which specifically concern the Navier-Stokes equations. In order to make the presentation shorter, we choose not to present  a number of results which have been proved by various authors under some additional geometrical assumptions on the flow, which imply the  conservation of quantities beyond the scaling (namely spherical, helicoidal or axisymmetric conditions). We refer for instance to \cite{ladyzhenskaya}, \cite{mahalovtitileibovich}, \cite{poncerackesideristiti}, or \cite{ukhoiudo} for such studies.

  To start with, let us recall the question asked in the previous paragraph, concerning the blow up of the~$\dot B^{-1+\frac dp}_{p,q}$ norm at blow-up time. A typical example of a solution with a finite~$\dot B^{-1+\frac dp}_{p,q}$ norm at blow-up time is a self-similar solution, and the question of the existence of such solutions was actually addressed by J. Leray in~\cite{leray}. The answer was given 60 years later by J. Ne\c{c}as, M. Ru\c{z}i\c{c}ka and V. \c{S}ver\'ak in\ccite{NRS}.  By analyzing the profile equation, they proved that there is no self-similar solution in~$L^3$ in three space dimensions. Later L. Escauriaza, G. Seregin and V. \c{S}ver\'ak were able to prove more generally that  if the solution is bounded in~$L^3 $, then it is regular (see~\cite{ESS}): in particular any solution blowing up in finite time must blow up in~$L^3 $.

Now let us turn to the existence of large, global unique solutions to the Navier-Stokes system in three space dimensions.
    
%
%We first notice that for axi-symmetric without swirl or helicoidal flows, there exist a unique global solution. This is a global existence result for large data but with some geometrical symmetries. Two heuristically reasons allow to understand why we have a global solutions in this case. First of all, the solution depends only on two variables and the two dimensional Navier-Stokes equations are globally well posed. Secondly, it is well known that the singularity set for weak solution is a set of Hausdorf measure smaller than $1/2$. If a symmetric flow has a singularity, then we have a curve of singularities which is of dimension Hausdorf 1, and in contradictions with theorem. In fact, This results used the conservation of some quantities which are more regular than the scaling regularity level.

An important example where  a unique global in  time solution exists for large initial data is the case where the domain   is thin in the vertical direction (in three space dimensions): that was proved by G. Raugel and G. Sell  in\ccite{raugelsell}
(see also the paper\ccite{iftimieraugelsell} by D. Iftimie, G. Raugel and G. Sell). 
Another example of large initial data generating a global solution was obtained by A. Mahalov and B. Nicolaenko in~\cite{babinandco}: in that case,  the initial
data is chosen so as to transform the equation into a rotating fluid equation (for which it is known that global solutions exist for a sufficiently strong rotation). 

In both those examples, the global wellposedness  of the two dimensional equation is an important ingredient in the proof. Two of the authors  also used such a property to construct in \cite{cgens} an example of periodic initial data which is large in $\dot B^{-1}_{\infty,\infty}$
but yet generates a global solution. Such an initial data is given by
$$
u_0^N(x) \eqdefa (Nu_h(x_h)\cos(Nx_3),-\dive_h u_h(x_h) \sin(Nx_3)),
$$
where $\|u_h\|_{L^2(\TT^2 )}\leq C(\ln N)^{\frac 14}$, and its $\dot B^{-1}_{\infty,\infty}$ norm is typically of the same size. This was generalized to the case of the space~$\R^3$ in~\cite{cg2}.

Similarly in\ccite{cg3},   that fact was used to prove a global existence result for large data which are slowly 
varying in one direction.   More precisely, if~$(v_0^h,0)$  and~$w_{0}$ are two smooth divergence free 
"profile" vector fields, then they proved that the initial data
\beq
\label{wellprepareddata}
 u_{0,\e} (x_h,x_3) \eqdefa (v_0^h(x_h,\e x_3), 0) + (\e w_0^h(x_h,\e x_3), w_0^3(x_h,\e x_3))
 \eeq
 generates, for~$\e$ small enough, a global smooth solution. Here, we
 have denoted~$x_h = (x_1,x_2)$. Using for instance the language of geometrical optics in the context
of fast rotating
 incompressible fluids,  and thinking of the problem in terms of   convergence to the two dimensional situation, 
this case can be seen as a "well prepared" case. We shall be coming back to that example in the next paragraph.

As a conclusion of this short (and of course incomplete) survey, let us present some results 
for the Navier-Stokes system with viscosity vanishing in the vertical direction. Analogous results to the classical Navier-Stokes system in the framework of small data are proved in\ccite{cdgg2},\ccite{dragosunique},\ccite{marius1} and\ccite{cz1}). To circumvent the difficulty linked with the absence of vertical viscosity, the key idea, which will be also crucial here (see for instance the proof of the second estimate of Proposition\refer{conventionalestimate}) is the following: 
the vertical derivative $\partial_3$  appears in the nonlinear term of the equation with the prefactor~$u_3$, which has some additional smoothness thanks to the divergence free condition which states that~$\partial_3 u_3 = - \partial_1 u_1 - \partial_2 u_2$.

 \subsection{Statement of the main result}
 In this work, we are interested in generalizing the situation\refeq{wellprepareddata} to the ill prepared  case: we shall investigate the case of initial data of the form 
 $$
  \Bigl(v_0^h(x_h,\e x_3) ,\frac 1 \e v_0^3(x_h,\e x_3)\Bigr),
 $$
 where~$x_h$ belongs to the torus~$\TT^2$ and~$x_3$ belongs to~$\R$. 
The main theorem of this article is the following.
\begin{theo}
\label{CGPtheo0}
{\sl Let~$a$ be a positive number. There are two positive numbers~$\e_0$ and~$\eta$  such that for any divergence free vector field~$v_0$ satisfying
$$
\|e^{a|D_3|} v_0\|_{H^{4}} \leq \eta,
$$
then, for any positive $\e$ smaller than~$\e_0$, the initial data
$$
u_{0,\e} (x)\eqdefa \Bigl(v_0^h(x_h,\e x_3) ,\frac 1 \e v_0^3(x_h,\e x_3)\Bigr)
$$
generates a global smooth solution of~$(NS)$ on~$\TT^2\times\R$.
}
\end{theo}

{\bf Remarks } 
\begin{itemize} 
\item
Such an initial data may be arbitrarily large in~$\dot B^{-1}_{\infty,\infty}$, more  precisely of size~$\e^{-1}$. Indeed it is proved in~\cite{cg3}, Proposition~1.1, that if~$f$ and~$g$ 
are two functions   in~$  \cS(\TT^2)$ and $\cS(\R)$ respectively, then~$h^\e (x_h,x_3) \eqdefa f(x_h)g(\e x_3)$ satisfies, if $\e$ is small enough,
$$
\|h^\e\|_{\dot B^{-1}_{\infty,\infty} } \geq \frac 1 4 \|f\|_{\dot B^{-1}_{\infty,\infty}} 
\|g\|_{L^\infty}.
$$

\item As in the well prepared case studied in~\cite{cg3} and recalled in the previous paragraph, the structure of the nonlinear term will have a crucial role to play in the proof of the theorem.

\item
The reason why the horizontal variable is restricted to a torus is to be able to deal with very low horizontal frequencies: as it will be clear in the proof of the theorem, functions with zero horizontal average are treated differently to the others, and it is important that no small horizontal frequencies appear other than zero. In that situation, we are able to solve globally in time the equation (conveniently rescaled in~$\e$) for small  analytic-type initial data. We recall that in that spirit, some local in time results for Euler and Prandtl equation with analytic initial data can be found in\ccite{sammartino&caflisch}.  In this paper we shall follow a method close to a method introduced in\ccite{chemin21}.
 
\item
We finally  note that we can add to our initial data any small enough  data in $\dot H^{\frac 12}$, and we still obtain the global existence of the solution. Indeed, by the results contained in\ccite{gallagheriftimieplanchon}, if we fix an initial data which gives a global in time solution, then, all initial data in a small neighborhood,  give global in time solutions.
 \end{itemize}
 
 \section*{Acknowledgments}

The authors wish to thank Vladimir   \c{S}ver\'ak for pointing out the interest of this problem to them. They also thank   Franck Sueur for suggesting the analogy with Prandlt's problem.

\section{Structure of the proof }\label{mainsteps}

\setcounter{equation}{0} 

%%%%%%%%%%%%%%%%%%%%%%%%%%%%%%%%%%%%%%%%%%%%%%%%%%%%%%%%%%%%%%%%%%%%%%%%%%%%%%%%%%%%%%%%%%%%%%%%%%%%%

 \subsection{Reduction to a rescaled problem}\label{reduction}
We look for  the solution under  the form
$$
u_{\e} (t,x)\eqdefa \Bigl(v^h(t,x_h,\e x_3) ,\frac 1 \e v^3(t,x_h,\e x_3)\Bigr).
$$
This leads to the following  rescaled Navier-Stokes system.
$$
(RNS_\e)\ \left\{
\begin{array}{c}
\partial_{t} v^h -\Delta_\e v^h +v\cdot\nabla v^h=-\nabla^h q\\
\partial_{t} v^3 -\Delta_\e v^3 +v\cdot\nabla v^3=-\e^2\partial_3 q\\
\dive v =0\\
v_{|t=0}=v_{0}
\end{array}
\right.
$$
with~$\D_\e\eqdefa \partial_1^2+\partial_2^2+\e^2\partial_3^2$. As there is no boundary, the rescaled pressure~$q$ can be computed with the formula
\beq
\label{eqrescaledpressure}
\D_\e q = \sum_{j,k} \partial_j v^k\partial_k v^j =  \sum_{j,k} \partial_j \partial_k  (v^jv^k).
\eeq
It turns out that when~$\e$ goes to~$0$,~$\D^{-1}_\e$ looks like~$\D_h^{-1}$. In the case of~$\R^3$, for low horizontal  frequencies,  an expression of the type~$\D_h^{-1}(ab)$  cannot be estimated in~$L^2$ in general. This is the reason why we work in~$\TT^2\times\R$. In this domain, the problem of low horizontal frequencies reduces to the problem of the horizontal average that we denote by
$$
(M f )(x_3)\eqdefa \overline  f(x_3) \eqdefa \int_{\TT^2} f(x_h,x_3) dx_h.
$$ 
Let us also define~$M^\perp f\eqdefa(\Id-M)f$.
Notice that, because the vector field~$v$  is divergence free, we have~$\vbar ^3\equiv 0$.
The system~$(RNS_\e)$ can be rewritten in the following form.
$$
(RNS_\e)\ \left\{
\begin{array}{c}
\ds \partial_{t} w^h -\Delta_\e w^h + M^\perp\bigl(v \cdot\nabla w^h+w^3\partial_3\vbar^h\big)=-\nabla^h q\\
\ds \partial_{t} w^3 -\Delta_\e w^3 + M^\perp (v\cdot\nabla w^3)=-\e^2\partial_3 M^\perp q\\
\ds \partial_{t} \vbar^h -\e^2\partial_3^2 \vbar^h =-\partial_3M( w^3w^h) \\
\dive (\vbar +w) =0\\
(\vbar, w)_{|t=0}=(\vbar_0, w_{0}).
\end{array}
\right.
$$
The problem to solve this sytem is that there is no obvious way to compensate the loss  of one vertical derivative which appears in the 
equation on~$w_h$ and~$\vbar$ and also, but more hidden, in the pressure term. The method we use is inspired by the one introduced in\ccite{chemin21}
 and can be understood as a global Cauchy-Kowalewski  result. This is the reason why the hypothesis of analyticity in the vertical variable is required in our theorem.
  
Let us denote by~$\cB$ the unit ball of~$\R^3$ and by~$\cC$ the annulus of small radius~$1$ and large radius~$2$. For non negative~$j$, let us denote by~$L^2_j$ the space~$\cF L^2((\ZZ^2\times \R)\cap 2^j\cC)$ and by~$L^2_{-1}$ the space $\cF L^2((\ZZ^2\times \R)\cap \cB)$ respectively equipped with the (semi) norms
$$
\|u\|^2_{L^2_j} \eqdefa (2\pi)^{-d} \int_{2^j\cC} |\wh u(\xi)|^2d\xi \andf 
\|u\|^2_{L^2_{-1}} \eqdefa (2\pi)^{-d} \int_{\cB} |\wh u(\xi)|^2d\xi.
$$  
Let us now recall the definition of
inhomogeneous Besov spaces modeled on~$L^2$. 
\begin{defi}
{\sl Let $s$ be a nonnegative real number. The space~$B^s$  is the subspace of~$L^2$ such that
$$ 
\| u\|_{ B^s}\eqdefa \bigl\| \bigl( 2^{js} \| u\|_{L_j^2}\bigr)_j\bigr\|_{\ell^1}<\infty.
$$ 
}
\end{defi} 
We note that $u\in B^s$ is equivalent to writing~$\|u\|_{L_j^2}\leq C c_j 2^{-js}\|u\|_{B^s}$
where $(c_j )$ is a non negative series which belongs to the sphere of~$\ell^1$. Let us notice that~$B^{\frac 3 2}$ is included in~$\cF(L^1)$ and thus in the space of continuous bounded functions. Moreover, if we substitute $\ell^2$ to~$\ell^1$ in the above definition,  we recover the classical Sobolev space~$H^s$. 
 
 The theorem we actually prove is the following.
\begin{theo}
\label{CGPtheo1}
{\sl Let~$a$ be  a positive number. There are two positive numbers~$\e_0$ and~$\eta$  such that for any divergence free vector field~$v_0$  satisfying
$$
\|e^{a|D_3|} v_0\|_{B^{\frac 7 2}} \leq \eta,
$$
then, for any positive $\e$ smaller than~$\e_0$, the initial data
$$
u_{0,\e} (x)\eqdefa \Bigl(v_0^h(x_h,\e x_3) ,\frac 1 \e v_0^3(x_h,\e x_3)\Bigr)
$$
generates a global smooth solution of~$(NS)$ on~$\TT^2\times\R$.
}
\end{theo}

  \subsection{Definition of the functional setting}
\subsubsection{Study of a model problem}\label{modelproblem}
  In order to motivate the functional setting and to give a flavour of the method used to prove the theorem, let us study for a moment the following simplified model problem for~$(RNS_\e)$, in which we shall see in a rather easy way how the same type of  method as that of\ccite{chemin21} can be used (as a global Cauchy-Kowaleswski technique): the idea is to control a nonlinear quantity, which depends on the solution itself. So let us consider the equation
 $$
\partial_t u+\gamma u+a(D) (u^2)=0,$$ 
where~$u$ is a scalar, real-valued function, $\gamma$ is a positive parameter,   and $a(D)$ is a Fourier multiplier
 of order one.
We shall sketch the proof of the fact   that if the initial data satisfies, for some positive~$\delta$ and some small enough constant~$c$,
$$
\|u_0\|_{X}\eqdefa
\int e^{\delta |\xi|} |\widehat u(\xi) | \: d\xi\leq c\gamma,
$$
 then one has a global smooth solution, say in the space~$ \mathcal F(L^1)$ as well as all its derivative.  The idea of the proof is the following: we want to control the same 
 kind of quantity on the solution, but one expects
 the radius of  analyticity of the solution to decay in time.   So let  us introduce $\theta(t)$ the "loss of analyticity" of the solution, solving the following ODE: 
$$
 \dot\theta (t)\eqdefa \int e^{(\delta- \lambda\theta(t))|\xi|}|\widehat u(\xi)|d\xi \quad \mbox{with} \quad \theta(0)=0.
$$
The parameter~$\lambda$ will be chosen large enough at the end
   and we will prove that $\delta - \lambda\theta (t)$ remains positive for all times. The computations that follow hold as long as that assumption is true (and a bootstrap will prove that in fact it does remain true for all times).
We define the notation
$$
u_\theta (t) = {\mathcal F}^{-1} \left(e^{(\delta - \lambda\theta (t))  |\cdot| }  |  \: \widehat u(t,\cdot)|\right).
$$
Notice that
\beq\label{defthetaexample}
  \dot\theta (t) = \|u_\theta (t) \|_{ \mathcal F(L^1)} \quad\text{and}\quad\theta(t)=\int_0^t\|u_\theta (t')\|_{\mathcal F(L^1)}\,dt'.
\eeq
Taking the Fourier transform of the equation gives
$$
|\widehat u(t,\xi)| \leq e^{-\gamma t} |\widehat u_0(\xi)|+C \int_0^t e^{-\gamma(t-t') }  \:  |\xi|  \:  | {\mathcal F} (u^2) (t',\xi) | \: dt'.
$$
Using the fact that 
$$
(\delta - \lambda\theta (t))  \:  |\xi| \leq (\delta -\lambda \theta (t'))   \:  |\xi - \eta| + (\delta -\lambda \theta (t'))   \:  |  \eta| - \lambda  |\xi|  \: \int_{t'}^{t} \dot\theta (t) \: dt''\,,
$$
we infer that
$$
|\widehat u_\theta (t,\xi)|  \leq e^{-\gamma t} e^{\delta |\xi|}|\widehat u_0(\xi)|+C \int_0^t e^{-\gamma(t-t') - \lambda  |\xi| \int_{t'}^{t} \dot\theta (t) \: dt''}    |\xi| \,
  | {\mathcal F} (u_\theta^2) | (t',\xi) \: dt'.
$$

We note the important fact that 
$$
\int_0^t e^{ - \lambda  |\xi| \int_{t'}^{t} \dot\theta (t) \: dt''}    |\xi| \, \dot \theta (t') \, dt'\leq \frac{C}{\lambda}\,\virgp
$$
 which is very useful in what follows.
As ${\mathcal F}(ab) = (2\pi)^{-d}(\wh a\star \wh b)$, we have, for any~$t'\leq t$,
$$
| {\mathcal F} (u_\theta^2) | (t',\xi) \leq
\Bigl( \sup\limits_{0\leq t'\leq t} |\widehat u_\theta (t',\cdot )|\Bigr)\star |\widehat
u_\theta (t',\cdot )|.
$$
Recalling that~${ \mathcal F}(L^1)$ is an algebra, we infer that
$$
\Big\|\sup\limits_{0\leq t'\leq t} |\widehat u_\theta (t')|\Big\|_{L^1}\leq
\|u_0\|_{X}+\frac{C}{\lambda}\Big\|\sup\limits_{0\leq t'\leq t} |\widehat
u_\theta (t')|\Big\|_{L^1}
$$
and
$$
\theta(t)\leq C_\gamma\Big(\|u_0\|_{X}+\Big\|\sup\limits_{0\leq t'\leq t}
|\widehat u_\theta (t')|\Big\|_{L^1}\theta(t)\Big).
$$ 
%  It is then an easy exercise to integrate in~$\xi$ and take $L^\infty$ and $L^1$ norms in time (recalling that~$ \mathcal F(L^1)$ is an algebra) to find
%$$
%\Big\|\sup\limits_{0\leq t'\leq t} |\widehat u_\theta (t')|\Big\|_{L^1}\leq \|u_0\|_{X}+\frac{C}{\lambda}\Big\|\sup\limits_{0\leq t'\leq t} |\widehat u_\theta (t')|\Big\|_{L^1}
%$$
%and 
%$$
%\|\dot \theta\|_{L^1_t}\leq C_\gamma\bigg(\|u_0\|_{X}+\frac{C}{\lambda}\Big\|\sup\limits_{0\leq t'\leq t} |\widehat u_\theta (t')|\Big\|_{L^1}\|\dot \theta\|_{L^1_t}\bigg).
%$$
This allows by bootstrap to obtain the global in time existence of the solution, as soon as the initial data is small enough; we skip the computations, as they will be presented in full detail for the case of the system~$(RNS_\e)$.

\subsubsection{Functional setting}
In the light of the computations of the previous section, let us introduce the functional setting we are going to work with to prove the theorem.   
The proof relies on exponential decay estimates for the Fourier transform of the solution. Thus, for any locally bounded function~$\Psi$ on~$\R^+\times \ZZ^2\times\R$ and for any  function~$f$,  continuous in time and compactly supported in Fourier space, we define
$$
(f_\Psi)(t) \eqdefa \cF^{-1} \bigl(e^{\Psi(t,\cdot)}\wh f(t, \cdot)\bigr).
$$

Now let us define the key quantity we wish to control in order to prove the theorem. In order to do so, let us consider the Friedrichs approximation of the original $(NS)$ system
$$
\begin{cases}
\partial_t u- \Delta u+\P_n (u\cdot\nabla u+\nabla p)=0 \\
\dive u=0\\
u|_{t=0}=\P_n u_{0,\e},
\end{cases}
$$
where~$\P_n$ denotes the orthogonal projection of~$L^2$ on functions the Fourier transform of which is supported in the ball~$B_n$ centered at the origin and of radius~$n$.  Thanks to the $L^2$ energy estimate, this approximated system has a global solution the Fourier transform of which is supported in~$B_n$. Of course, this 
provides an approximation of the rescaled system namely 
$$
(RNS_{\e,n})\ \left\{
\begin{array}{c}
\ds \partial_{t} w^h -\Delta_\e w^h + \P_{n,\e} M^\perp\bigl(v\cdot\nabla w^h+w^3\partial_3\vbar+\nabla^h q\bigr)=0\\
\ds \partial_{t} w^3 -\Delta_\e w^3 + \P_{n,\e} M^\perp\bigl (v\cdot\nabla w^3+\e^2
\partial_3  q\bigr)=0\\
\ds \partial_{t} \vbar^h -\e^2\partial_3^2 \vbar^h +\P_{n,\e} \partial_3M( w^3w^h)=0 \\
\dive (\vbar +w) =0\\
(\vbar, w)_{|t=0}=(\vbar_0, w_{0}),
\end{array}
\right.
$$
where   $\P_{n,\e}$ denotes the orthogonal projection of~$L^2$ on functions the Fourier transform of which is supported in~$B_{n,\e} \eqdefa \{ \xi\,/\ |\xi_\e|^2\eqdefa |\xi_h|^2+\e^{2} \xi_3^2 \leq n^2\}$. We shall prove analytic type estimates here, meaning exponential decay estimates for the the solution of the above approximated system. In order to make notation not too  heavy we will drop the fact that the solutions we deal with are in fact approximate solutions and not solutions of the original system.
A priori bounds on the approximate sequence will be derived, which will clearly yield the same bounds on the solution. In the spirit of\ccite{chemin21} (see also\refeq{defthetaexample} in the previous section), we define the function~$\theta$  (we drop also the fact that~$\theta$  depends on~$\e$ in all that follows) by
\beq
\label{definfunctionfund}
\dot \theta (t) = \|w^3_\Phi(t)\|_{B^{\frac 7 2}} +\e \|w^h_\Phi(t)\|_{B^{\frac 7 2}}\andf
\theta(0)=0
\eeq
where
\beq
\label{definfunctionfundeq1}
\Phi(t,\xi) = t^{\frac 1 2} |\xi_h| +a|\xi_3| -\lam \theta (t) |\xi_3| 
\eeq
for some~$\lam$ that will be chosen later on (see Section\refer{prooftheorem}). Since the Fourier transform of~$w$ is compactly supported, the 
above differential equation  has a unique global solution on~$\R^+$.  If we prove that 
\beq
\label{ineganalyticfund}
\forall t\in \R^+\,,\ \theta (t) \leq \frac  a \lam\,\virgp
\eeq
this will imply that the sequence of approximated solutions of the rescaled system is a bounded sequence of~$L^1(\R^+; \mbox{Lip})$. So is, for a fixed $\e$, the family of approximation  of the original Navier-Stokes equations. This is   (more than) enough to imply that a global smooth solution exists.

%%%%%%%%%%%%%%%%%%%%%%%%%%%%%%%%%%%%%%%%%%%%%%%%%%%%%%%%%%%%%% 
\subsection{Main steps of the proof}
The proof of Inequality\refeq{ineganalyticfund} will be a consequence of the  following  two propositions  which provide estimates on~$v^h$, $w^h$ and~$w^3$. For technical reasons, these statements require the  use  of a modified version (introduced in\ccite{chemin13}) of~$L^\infty_T(B^s)$ spaces.
\begin{defi}
\label{definLinftytilde}
{\sl Let~$s$ be a real number.  We define the space~$\wt L^\infty_T(B^s)$ as the subspace  of   functions~$f$ of~$ L^\infty_T(B^s)$ such that the following quantity is finite:
$$
\|f\|_{\wt L^\infty_T(B^s)} \eqdefa \sum_j 2^{js} \|f\|_{L^\infty_T(L_j^2)}.
$$}
\end{defi}

Theorem\refer{CGPtheo1} will be an easy consequence of the following propositions, which will be proved in the coming sections.

The first one  uses only the fact that the function $\Phi$ is subadditive. 

\begin{prop}
\label{conventionalestimate}
{\sl A constant~$C_0^{(1)}$ exists such that, for any positive~$\lam$,   for any initial data~$v_0$, and 
 for any~$T$ satisfying~$\displaystyle  \theta (T) \leq  a / \lam $, we have
$$
\theta(T) \leq \e \|e^{a|D_3|}w^h_0\|_{B^{\frac 7 2}}+ \|e^{a|D_3|}w^3_0\|_{B^{\frac 7 2}} + C_0^{(1)} \|v_\Phi\|_{\wt L^\infty_T(B^\frac 7 2)}\displaystyle  \theta (T) .
$$
Moreover, we have the following~$L^\infty$-type estimate on the vertical component:
$$
\|w^3_\Phi\|_{\wt L^\infty_T(B^{\frac 7 2})}  \leq \|e^{a|D_3|}w^3_0\|_{B^{\frac 7 2}} + C_0^{(1)} \|v_\Phi\|_{\wt L^\infty_T(B^{\frac 7 2})}^2 .
$$
}
\end{prop}
The second one is   more subtle to prove, and it shows that the use of the analytic-type norm
actually allows to recover the missing vertical derivative on~$v^{h}$, in a~$L^{\infty}$-type space. It 
should be compared to the methods described in Section~\ref{modelproblem}.

\begin{prop}
\label{estimhorinzonLinftytilde}
{\sl A constant~$C_0^{(2)}$  exists such that, for any positive~$\lam$,  for any initial data~$v_0$, and 
 for any~$T$ satisfying~$\displaystyle   \theta (T) \leq  a / \lam $, we have
$$
\|v_\Phi^h\|_{\wt L^\infty_T(B^{\frac 7 2})} \leq  \|e^{a|D_3|} v^h_0\|_{B^{\frac 7 2}} + C_0^{(2)} \Bigl( \frac 1 \lam 
+ \|v_\Phi\|_{\wt L^\infty_T(B^{\frac 7 2})}\Bigr)\|v_\Phi^h\|_{\wt L^\infty_T(B^{\frac 7 2})} .
$$
}
\end{prop}

%%%%%%%%%%%%%%%%%%%%%%%%%%%%%%%%%%%%%%%%%%%%%%%%%%%%%%%%%%%%%% 

\subsection{Proof of the theorem assuming the two propositions}
\label{prooftheorem}
Let us assume these two propositions are true 
for the time being and conclude the proof of Theorem\refer{CGPtheo1}. It relies on a continuation argument.

For any positive~$\lam$ and~$\eta$, let us define
$$
 \cT_\lam\eqdefa \bigl\{T\,/\ \max \{ \|v_\Phi\|_{\wt L^\infty_T(B^\frac72)},\theta (T)\} \leq 4\eta\bigr\},
$$
 As the two functions involved in the definition of~$\cT_\lam$ are non decreasing,~$\cT_\lam$ is an interval. 
As~$\theta$ is an increasing function which vanishes at~$0$, a positive time~$T_0$ exists
 such that~$\theta(T_{0})\leq 4\eta$. 
 Moreover, if~$\|e^{a|D|_3|}v_0\|_{B^{\frac 7 2}} \leq \eta$  then, since~$\partial_tv =\P_n F(v)$ (recall 
that we are considering Friedrich's approximations), a positive time~$T_1$ (possibly 
depending on~$n$) exists such that~$ \|v_\Phi\|_{\wt L^\infty_{T_{1}}(B^{\frac 7 2})}\leq 4\eta$. Thus~$\cT_\lam$ is  the form~$[0,T^\star)$ for some positive~$T^\star $. Our purpose is to prove that~$T^\star=\infty$. As we want to apply Propositions~\ref {conventionalestimate} and~\ref{estimhorinzonLinftytilde}, we need that~$\lam \theta (T)\leq a$. This leads to the condition
 \beq
 \label{condPhipositive}
 4\lam\eta \leq a .
 \eeq
From Proposition\refer{conventionalestimate}, defining~$C_0\eqdefa C_0^{(1)}+C_0^{(2)}$, we have, for all~$T\in \cT_\lam$,
$$
\|v_\Phi\|_{\wt L^\infty_T(B^{\frac 7 2})}  \leq \|e^{a|D_3|} v_0\|_{B^{\frac 7 2}} +
\frac {C_0} \lam \|v_\Phi\|_{\wt L^\infty_T(B^{\frac 7 2})} + C_0\|v_\Phi\|^2_{\wt L^\infty_T(B^{\frac 7 2})}.
$$
Let us choose~$\lam= \ds \frac 1 {2C_0}\,\cdotp$ This gives
$$
\|v_\Phi\|_{\wt L^\infty_T(B^{\frac 7 2})}  \leq 2\|e^{a|D_3|} v_0\|_{B^{\frac 7 2}} +
4 C_0\eta \|v_\Phi\|_{\wt L^\infty_T(B^{\frac 7 2})}.
$$
Choosing $\ds \eta = \frac 1 {12 C_0}\virgp$ we infer that, for any~$T\in \cT_\lam$,
\beq
\label{proofconcludeeq1}
\|v_\Phi\|_{\wt L^\infty_T(B^{\frac 7 2})}  \leq 3\|e^{a|D_3|} v_0\|_{B^{\frac 7 2}} .
\eeq
 Propositions\refer{conventionalestimate} and~\ref{estimhorinzonLinftytilde} imply that, for all $T\in \cT_\lam$,
$$
\theta(T) \leq \e \|e^{a|D_3|}w^h_0\|_{B^{\frac 7 2}}+ \|e^{a|D_3|}w^3_0\|_{B^{\frac 7 2}} + C_0\eta \theta (T) .
$$
This implies that
$$
\theta(T) \leq 2\e \|e^{a|D_3|}w^h_0\|_{B^{\frac 7 2}}+ 2\|e^{a|D_3|}w^3_0\|_{B^{\frac 7 2}}.
$$
If~$2\e \|e^{a|D_3|}w^h_0\|_{B^{\frac 7 2}}+ 2\|e^{a|D_3|}w^3_0\|_{B^{\frac 7 2}} \leq \eta$ and
$\|e^{a|D_3|}v^h_0\|_{B^{\frac 7 2}} \leq \eta$, then the above estimate and Inequality\refeq{proofconcludeeq1} ensure\refeq{ineganalyticfund}. This concludes the proof of Theorem\refer{CGPtheo1}.

%%%%%%%%%%%%% NOTATION AND PRELIMINARY RESULTS %%%%%%%%%%% 

\section{The action of subadditive phases on (para)products}
\label{preliminaryresults}
\setcounter{equation}{0} 

%%%%%%%%%%%%%%%%%%%%%%%%%%%%%%%%%%%%%%%%%%%%%%%%%%%

It will be useful to consider, for any function~$f$,  the inverse Fourier transform of~$|\wh f|$, defined~as
$$
f^+\eqdefa \cF^{-1}|\wh f|.
$$
 Let us notice that the map~$f\mapsto f^+$ preserves the norm of all~$B^s$ spaces.
In all this section, $\Psi$ will denote a locally bounded function on~$\R^+\times \TT^2\times\R$ which satisfies the following inequality
\beq
\label{subadditivitycond}
\Psi(t,\xi) \leq \Psi(t,\xi-\eta)+\Psi(t,\eta).
\eeq

In all the following, we will denote by~$C$ or~$c$ universal constants, which do not depend on any of the parameters of the problem, and which
may change from line to line. We will denote generically by~$c_{j}$ any sequence in~$\ell^{1}(\ZZ)$ of norm~$1$.
 
We shall denote by~$\EE f$   
the solution of~$\partial_t g-\D_\e g= f$   
with initial data equal to~$0$. 
 We use also a very basic version of  Bony's decomposition. Let us define (using the notation introduced in Section~\ref{reduction}),
 $$
 T_ab \eqdefa \cF^{-1} \sum_{j}\int_{2^{j}\cC \cap \cB(\xi, 2^{j}) } \!\!\wh a(\xi-\eta) \wh b(\eta) d\eta \!\!\andf \!\! R_{a}b \eqdefa  \cF^{-1} \sum_{j}\int_{2^{j}\cC \cap \cB(\xi, 2^{j+1}) } \!\!
 \wh a(\xi-\eta) \wh b(\eta) d\eta .
$$
We obviously have   $ab = T_ab+R_{b}a$. 

The way  the Fourier multiplier~$e^\Psi$ acts on bilinear functionals is described by the following lemma.
\begin{lemma}
\label{estmitildepastilde}
{\sl For any positive~$s$, a constant~$C$ exists which satisfies the following properties. For any function~$\Psi$ satisfying\refeq{subadditivitycond}, for any function~$b$ in~$  L^1_T(B^{s})$, a positive sequence~$\suite c j \ZZ$ exists in the sphere of~$\ell^1(\ZZ)$ such that, for any~$a$ in~$   L^1_T(B^{\frac 3 2})$, and any~$t\in [0,T]$, we have
$$
\|(T_ a b )_\Psi (t) \|_{L_j^2} + \|( R_ a b)_\Psi  (t) \|_{L_j^2} \leq  C
c_j 2^{-js} \|a(t)\|_{B^{\frac 3 2}}  \min\bigl \{ \|b(t)\|_{B^{s}} , \|b\|_{\wt L^\infty_T(B^s)}\bigr\}.
$$
}
\end{lemma}

\begin{proof} 
We prove only the lemma for~$R$, the proof for~$T$ being strictly identical. Let us first investigate the
 case when the function~$\Psi$ is identically~$0$.  We first observe that for any~$\xi$ in the annulus~$2^j\cC$, we have
$$
\cF (R_a b(t))(\xi) = \sum_{j'\geq j-2}\int_{2^{j'}\cC \cap B(\xi, 2^{j'+1}) } \wh a(t,\xi-\eta) \wh b(t,\eta) d\eta .
$$
As~$B^{\frac 3 2}$ is included in~$\cF(L^1)$, we infer that, by definition of~$\|\cdot\|_{\wt L^\infty_T(B^s)}$,
$$
\|R_a b(t)\|_{L_j^2}\leq C \|a(t)\|_{B^{\frac 3 2}}\sum_{j'\geq j-2}c_{j'}2^{-j's} \min \bigl\{ \|b(t)\|_{B^{s}} , \|b\|_{\wt L^\infty_T(B^s)}\bigr\}.
$$
Defining~$\ds \widetilde c_j=\sum\limits_{j'\geq j-2}\ 2^{(j-j')s}c_{j'}$ which satisfies $\ds\sum_j \tilde c_j\leq C_s$, we obtain
\beq
\label{estmitildepastildedemoeq1}
\| R_a b(t)\|_{L_j^2}\leq C \tilde c_j2^{-js}\|a(t)\|_{B^{\frac 32}} \min \bigl\{ \|b(t)\|_{B^{s}} , \|b\|_{\wt L^\infty_T(B^s)}\bigr\}.
\eeq
The lemma is then proved in the case when the function~$\Psi$ is identically~$0$.  In order to treat the general case, let us write that
\begin{eqnarray*}
 e^{\Psi(t,\xi)}\cF(R_a b) (\xi) & = &   e^{\Psi(t,\xi)} \sum_{j} \int_{2^j\cC \cap B(\xi, 2^j) } \wh a(\xi-\eta) \wh b(\eta) d\eta\\
 & \leq &  \sum_{j} \int_{2^j\cC \cap B(\xi, 2^j)}  e^{\Psi(t,\xi-\eta)} \wh a^+(\xi-\eta) e^{\Psi(t,\eta)}\wh b^+(\eta)d\eta.
\end{eqnarray*}
This means exactly that~$|\cF(R_a b)_\Psi (t,\xi)| \leq \cF (R_{a^+_\Psi} b^+_\Psi)(t,\xi)$. Then, the estimate\refeq{estmitildepastildedemoeq1} implies the lemma.
\end{proof}

\begin{corol}
\label{paradiffanalytic}
{\sl If~$s$ is positive, we have, for any function~$\Psi$ satisfying\refeq{subadditivitycond},
 \begin{eqnarray*}
\|(T_a b)_\Psi\|_{\wt L^\infty_T(B^s)} + \|(R_a b)_\Psi\|_{\wt L^\infty_T(B^s)} & \leq
& C \|a_\Psi\|_{L^\infty_T(B^{\frac 32})} \|b_\Psi\|_{\wt L^\infty_T(B^s)} \andf \\ 
\|(T_a b)_\Psi\|_{L^1_T(B^s)} + \|(R_a b)_\Psi\|_{L^1_T(B^s)} & \leq
& C  \min\bigl\{ \|a_\Psi\|_{L^1_T(B^{\frac 32})} \|b_\Psi\|_{ \wt L^\infty_T(B^s)}, \\
&& \quad \quad \quad \quad \quad \quad 
\|a_\Psi\|_{L^\infty_T(B^{\frac 32})} \|b_\Psi\|_{ L^1_T(B^s)}\bigr\}.
\end{eqnarray*}
}
\end{corol}
\begin{proof} Taking the~$L^\infty$ norm in time on the inequality of Lemma\refer{estmitildepastilde} gives that
$$
\|(T_ a b )_\Psi \|_{L^\infty_T(L_j^2)} + \|( R_ a b)_\Psi \|_{L^\infty_T(L_j^2)} \leq  C
c_j 2^{-js} \|a\|_{L^\infty_T(B^{\frac 3 2})} \|b\|_{\wt L^\infty_T(B^s)}.
$$
which is the first inequality of the corollary. The proof of the second one is analogous.
\end{proof}
 
 %%%%%%%%%%%%%%%%%%%%% %%%%ACTION HEAT OPERATOR%%%%%%%%%%%%%%%%%%%%%%% 

\section{The action of   the phase~$\Phi$ on  the heat operator}
\setcounter{equation}{0} 

The purpose of this section  is the study of the action of the multiplier~$e^{\Phi}$ on $\EE f$. Let us recall that the function~$\Phi$ is defined in\refeq{definfunctionfundeq1} by
 ~$\Phi(t,\xi) = t^{\frac 1 2} |\xi_h| +a|\xi_3| -\lam \theta (t) |\xi_3| $. This action is described by the following lemma.
\begin{lemma}
\label{PhiactonEE}
{\sl A constant~$C_0$ exists such that, for any function~$f$ with compact spectrum, we have, for any~$s$,
\begin{eqnarray*}
 \nonumber \|(\EE M^\perp f)_\Phi\|_{\wt L^\infty_T(B^s)} & \leq & C_0  \|g_\Phi\|_{\wt L^\infty_T(B^s)}\andf\\
  \|(\EE M^\perp f)_\Phi\|_{L^1_T(B^s)} & \leq & C_0  \|g_\Phi\|_{L^1_T(B^s)} 
\quad\hbox{ where }\quad
g  \eqdefa  \cF^{-1}\Bigl( \frac 1 {|\xi_h|} |\cF M^\perp f|\Bigr).
 \end{eqnarray*} }
\end{lemma}
\begin{proof} 
Let us write~$\EE$ in terms of the Fourier transform. We have, for any~$\xi \in (\ZZ^2\setminus\{0\})\times\R$,
$$
\cF\left(\EE  f\right)_\Phi(t,\xi) = e^{\Phi(t,\xi)} \int_0^te^{-(t-t')|\xi_\e|^2} f(t',\xi)dt',
$$
with, as in all that follows,~$|\xi_\e|^2 \eqdefa |\xi_h|^2 + \e^2 |\xi_3|^2$.
Thus we infer, for any~$\xi \in (\ZZ^2\setminus\{0\})\times\R$,
$$
|\cF\left((\EE f)_\Phi\right)(t,\xi)| \leq  \int_0^t e^{-(t-t')|\xi_\e|^2+\Phi(t,\xi)-\Phi(t',\xi)} \cF(f_\Phi^+)(t',\xi)dt'.
$$
By definition  of $\Phi$, we have (see\ccite{chemin21}, estimate (24))
\beq
\label{relationLerayNantes}
\Phi(t,\xi)-\Phi(t',\xi) \leq -\lam |\xi_3| \int_{t'}^t \dot \theta(t'')dt'' +\frac {t-t'} 2 |\xi_h|^2.
\eeq
Thus we have, for any~$\xi \in (\ZZ^2\setminus\{0\})\times\R$,
\beq
\label{relationLerayNantescor1}
|\cF\left((\EE f)_\Phi\right)(t,\xi)| \leq  \int_0^t e^{-\frac {(t-t')} 2 |\xi_h  |^2 -\e^2(t-t')|\xi_3  |^2} \cF(f_\Phi^+)(t',\xi)dt'.
\eeq
Let us define~$\cC_h\eqdefa \{1\leq |\xi_h|\leq 2 \}\times \R$. The above inequality means that we have, for any~$\xi$ in~$2^j\cC\cap 2^k\cC_h$, $$
|\cF ((\EE f)_\Phi)(t,\xi)| \leq C \int_0^t  e^{-c(t-t')2^{2k}} 2^k
\wh g_\Phi(t',\xi) dt'. 
$$
Taking the~$L^2$ norm in~$\xi$ in that inequality gives
\beq
\label{PhiactonEEdemoeq1}
%\|(\EE f)_\Phi(t)\|_{\cF L^2(2^j\cC\cap 2^k\cC_h)} \leq \int_0^t  e^{-c(t-t')(2^{2k}+\e^22^{2j})} 2^k
%\|g_\Phi(t')\|_{L_j^2} dt'. 
\|(\EE f)_\Phi(t)\|_{\cF L^2(2^j\cC\cap 2^k\cC_h)} \leq \int_0^t  e^{-c(t-t') 2^{2k} } 2^k
\|g_\Phi(t')\|_{L_j^2} dt'. 
\eeq
By definition of the~$\wt L^\infty_T(B^s)$ norm, this gives, for any~$t\leq T$,
\begin{eqnarray*}
2^{js} \|(\EE f)_\Phi\|_{L^\infty_T(L^2(2^j\cC\cap 2^k\cC_h))} & \leq & C c_j \|g_\Phi\|_{\wt L^\infty_T(B^s)}  \int_0^t  e^{-c(t-t')2^{2k}} 2^k dt' \\
 & \leq  & C c_j2^{-k} \|g_\Phi\|_{\wt L^\infty_T(B^s)}.
\end{eqnarray*}
Now, writing that
$$
\|(\EE M^\perp f)_\Phi\|_{L^\infty_T(L_j^2)} \leq \sum_{k=0}^\infty  \|\EE(f_\Phi)\|_{L^\infty_T(L^2(2^j\cC\cap 2^k\cC_h))} 
$$
gives the first inequality of the lemma. 

In order to prove the second one, let us use the definition of the norm of the space~$B^s$ and~(\ref{PhiactonEEdemoeq1}); this gives
\begin{eqnarray*}
\sum_j 2^{js} \|(\EE f)_\Phi\|_{L^1_T(L_j^2)}  & \leq  & 
 \sum_{j,k}  2^{js} \|\EE(f_\Phi)\|_{L^1_T(\cF L^2(2^j\cC\cap 2^k\cC_h))} \\
& \leq  & C \sum_{j,k} \int_{[0,T]^2} {\mathbf 1}_{t\geq t'} e^{-c2^{2k}(t-t')} 2^k
c_j(t')\|g_\Phi(t')\|_{B^s} dt'dt.
\end{eqnarray*}
Integrating first in~$t$ gives
$$
\sum_j 2^{js} \|(\EE f_\Phi)\|_{L^1_T(L_j^2)}   \leq C \sum_{j,k} \int_{[0,T]} 2^{-k}
c_j(t')\|g_\Phi(t')\|_{B^s} dt.
$$
As the index~$k$  is nonnegative, we get the second estimate of the lemma. \end{proof}
The following lemma is a key one. It is here that the function~$\theta$ allows the gain  of the vertical derivative, in the spirit of the example presented in Section~\ref{modelproblem}. 
\begin{lemma}
\label{PhiactonEEd3gainlambda}
{\sl Let~$a(D)$  and $b(D)$ be   two Fourier  multipliers such that~$|a(\xi)|\leq C |\xi_3|$ and~$|b(\xi)|\leq C|\xi|^2$. We have
$$
\longformule{
\|(\EE a(D) R_{b(D)w^3} f)_\Phi\|_{\wt L^\infty_T(B^{\frac 7 2})}
+\|(\EE a(D) T_{b(D)w^3} f)_\Phi\|_{\wt L^\infty_T(B^{\frac 7 2})} 
}
{ 
{}
\leq
 C \Bigl(\frac 1 \lam +\|w^3_\Phi\|_{\wt L^\infty_T(B^{\frac 7 2})}\Bigr)
\|f_\Phi\|_{\wt L^\infty_T(B^{\frac 7 2})}.
}
$$ 
}
\end{lemma}
\begin{proof}
We give only the proof for the first term, the second term is estimated exactly along the same lines. Let us write $\EE$ in Fourier variables. We have
$$
\cF (\EE a(D)R_{b(D)w^3} f)_\Phi (t,\xi)=e^{\Phi(t,\xi)}\int_0^t e^{-(t-t')|\xi_\e|^2}a(\xi)\cF (R_{b(D)w^3}f)(t',\xi)dt'.
$$
Thus, using that $|a(\xi)|\leq C|\xi_3|$, we obtain
$$
|\cF(\EE a(D) R_{w^3}f)_\Phi(t,\xi)|\leq C\int_0^t e^{-(t-t')|\xi_\e|^2+\Phi(t,\xi)-\Phi(t',\xi)}|\xi_3|\,|\cF ((R_{b(D) w^3} f)_\Phi)(t',\xi)|dt'.
$$
Taking into account Inequality\refeq{relationLerayNantes}, we have
$$
|\cF(\EE a(D) R_{w^3}f)_\Phi(t,\xi)|\leq C\int_0^t e^{-\frac {t-t'} 2 |\xi_\e|^2-\lam|\xi_3|\int_{t'}^t\dot\theta (t'')dt''}|\xi_3|\,|\cF((R_{b(D)w^3} f)_\Phi)(t',\xi)|dt'.
$$
Let us denote by~$\Psi$ the Fourier multiplier~$\Psi a \eqdefa \cF^{-1}( {\bf 1}_{|\xi_h|\leq 2 |\xi_3|}\wh a)$. If~$|\xi_h|\leq 2 |\xi_3|$ and~$\xi$ is in~$2^j\cC$, then, we have that~$|\xi_3|\sim 2^j$. Thus, we infer that, for any~$\xi$ in~$2^j\cC$, 
$$
|\cF\Psi (\EE a(D) R_{b(D)w^3}f)_\Phi(t,\xi)|\leq \int_0^t e^{-c\lam2^j\int_{t'}^t\dot\theta (t'')dt''}
2^j |\cF((R_{b(D)w^3} f)_\Phi)(t',\xi)|dt'.
$$
Taking the $L^2$ norm gives
$$
\|\Psi (\EE a(D) R_{b(D)w^3}f)_\Phi(t,\cdot)\|_{L_j^2} \leq \int_0^t e^{-c\lam2^j\int_{t'}^t\dot\theta (t'')dt''}
2^j \|(R_{b(D)w^3} f)_\Phi)(t')\|_{L_j^2}dt'.
$$
Using Lemma\refer{estmitildepastilde}, we get
\begin{eqnarray*}
2^{j\frac 7 2} \|\Psi (\EE a(D) R_{b(D)w^3}f)_\Phi(t,\cdot)\|_{L_j^2} & \leq & 
C c_j \|f_\Phi(t)\|_{\wt L^\infty_T(B^{\frac 7 2})}  \\
&& \qquad{}\times\int_0^t e^{-c\lam2^j\int_{t'}^t\dot\theta (t'')dt''}
2^j \|b(D)w^3_\Phi(t',\cdot)\|_{B^{\frac 3 2}}dt' \\
& \leq & 
C c_j \|f_\Phi(t)\|_{\wt L^\infty_T(B^{\frac 7 2})}  \\
&& \qquad\quad{}\times \int_0^t e^{-c\lam2^j\int_{t'}^t\dot\theta (t'')dt''}
2^j \|w^3_\Phi(t',\cdot)\|_{B^{\frac 7 2}}dt' \\
& \leq & C c_j \|f_\Phi(t)\|_{\wt L^\infty_T(B^{\frac 7 2})}  \int_0^t e^{-c\lam2^j\int_{t'}^t\dot\theta (t'')dt''} 2^j \dot \theta(t') dt'\\
& \leq & \frac  C \lam c_j \|f_\Phi(t)\|_{\wt L^\infty_T(B^{\frac 7 2})}  .
\end{eqnarray*}
By summation in~$j$, we deduce that
\beq
\label{PhiactonEEd3gainlambdadeoeq1}
\|\Psi (\EE a(D) R_{w^3}f)_\Phi\|_{\wt L^\infty_T(B^{\frac 7 2})} \leq \frac  C \lam \|f_\Phi(t)\|_{\wt L^\infty_T(B^{\frac 7 2})}.
\eeq
If~$2|\xi_3|\leq |\xi_h|$, then, for any~$\xi$ in~$2^j \cC$, $|\xi_h|$ is equivalent to~$2^j$ 
and~$|\xi_{3}|$
is less than~$2^{j}$. So we infer that for any~$\xi$ in~$2^j\cC$, 
$$
|\cF(\Id-\Psi) (\EE a(D) R_{b(D)w^3}f)_\Phi(t,\xi)|\leq \int_0^t e^{-c(t-t')2^{2j} }
2^j |\cF((R_{b(D)w^3} f)_\Phi)(t',\xi)|dt'.
$$
By definition of~$\|\cdot\|_{\wt L^\infty_T(B^{\frac 7 2})}$, taking the $L^2$ norm of the above inequality gives 
\begin{eqnarray*}
2^{j\frac 7 2} \|(\Id-\Psi )(\EE a(D) R_{b(D)w^3}f)_\Phi(t,\cdot)\|_{L_j^2} & \leq & 
 \int_0^t e^{-c2^{2j}(t-t')} 2^j 2^{j\frac 7 2 }\|(R_{b(D)w^3} f)_\Phi)(t')\|_{L_j^2}dt'  \\ 
& \leq &  C c_j \|(R_{b(D)w^3} f)_\Phi\|_{\wt L^\infty_T(B^{\frac 7 2})}  .
\end{eqnarray*}
After a summation in~$j$, Corollary\refer {paradiffanalytic} implies that
\beno
 \|(\Id-\Psi )(\EE a(D) R_{b(D)w^3}f)_\Phi\|_{\wt L^\infty_T(B^{\frac 7 2})}  
 & \leq &
 C\|b(D)w^3_\Phi\|_{\wt L^\infty_T(B^\frac 3 2)} \|f_\Phi\|_{\wt L^\infty_T(B^{\frac 7 2})}\\
 & \leq &
 C\|w^3_\Phi\|_{\wt L^\infty_T(B^\frac 7 2)} \|f_\Phi\|_{\wt L^\infty_T(B^{\frac 7 2})} .
\eeno
Together with\refeq{PhiactonEEd3gainlambdadeoeq1}, this concludes the proof of the lemma.
\end{proof}

\begin{lemma}
\label{PhiactonEEepsd3}
{\sl A constant~$C_0$ exists such that, for any function~$f$ with compact spectrum, we have
 for~$\al$ in~$\{1,2\}$,
\begin{eqnarray*}
\bigl \|\bigl(\EE (\e \partial_3)^\al M^{\perp}f\bigr)_\Phi\bigr\|_{\wt L^\infty_T(B^s)} 
 &\leq & C_0  \|f_\Phi\|_{\wt L^\infty_T(B^s)}\andf \\
\bigl \|\bigl(\EE (\e \partial_3)^\al  M^{\perp} f\bigr)_\Phi\bigr\|_{L^1_T(B^s)}& \leq & C_0  \|f_\Phi\|_{L^1_T(B^s)}.
\end{eqnarray*}
 }
\end{lemma}
\begin{proof}
Let us start with the case when~$\al = 2$. Recalling\refeq{relationLerayNantescor1}, we have (for~$0<\e<1$),
$$
\e^2|\cF(\EE \partial_3^2 f)_\Phi(t,\xi)| \leq  \int_0^t e^{-\e^2\frac{(t-t')}2 |\xi|^2 } \e^2 \xi_3^2 \cF( f_\Phi^+)(t',\xi)dt'.
$$
Writing that $ |\xi_3| \leq |\xi|$, we infer that
$$
\e^2\bigl\| (\EE \partial_3^2 f)_\Phi (t)\bigr\|_{L_j^2}  \leq
\int_0^t e^{- c{\e^2} (t-t')2^{2j}} \e^2 2^{2j} \|f(t')\|_{L_j^2} dt'.
$$
The estimates follow directly by applying Young's inequality in $t$.

In the case when~$\al = 1$, we decompose~$ f$ into two parts, 
$$
f =  f^{(1)} + f^{(2)},\quad \mbox{with}\quad   f^{(1)} = \cF^{-1} ({\mathbf 1}_{\e |\xi_3| \leq |\xi_h|} \wh  f).
$$
Let us start by studying the first contribution. We simply write that
\begin{eqnarray*}
\e \bigl | \cF\bigl(\EE \partial_3  f^{(1)}\bigr)_\Phi(t,\xi)\bigr| &\leq& 
 \int_0^t e^{-\frac {(t-t')} 2 |\xi_\e|^2 } \e |\xi_3| \cF(  f^{(1)}_\Phi)^+(t',\xi)dt' \\
 &\leq& 
 \int_0^t e^{-\frac {(t-t')} 2 |\xi_\e|^2 }  |\xi_h| \cF(  f^{(1)}_\Phi)^+(t',\xi)dt' 
\end{eqnarray*}
which amounts exactly to the computation~(\ref{PhiactonEEdemoeq1}), with~$g$ replaced by~$  f^{(1)}$.  
On the other hand, for~$  f^{(2)}$ we can write  
$$
\widehat g^{(2)} (\xi) \eqdefa  \frac 1 {|\xi_h|}  {\mathbf 1}_{\e |\xi_3| \geq |\xi_h|}
| \cF M^{\perp}f^{(2)}(\xi)| 
$$
so that
$$
\e \Bigl |\cF \bigr(\EE \partial_3  M^{\perp}f^{(2)}\bigr)_\Phi (t,\xi)\Bigr | \leq \int_0^t e^{-\frac {(t-t')} 2 |\xi_h|^2 -\e^2(t-t') |\xi_3|^2 }
 \e |\xi_3| |\xi_h|   \widehat g^{(2)} (t',\xi)
dt' .
$$
Since~$|\xi_h| \leq \e |\xi_3| $, we are reduced to the case when~$\alpha = 2$ and the conclusion comes from the fact that~$\|g^{(2)}_\Phi\|_{B^s} \leq
\| M^{\perp}f^{(2)}_\Phi\|_{B^s} \leq \| f_\Phi\|_{B^s} $.
That proves the lemma.
\end{proof}

%%%%%%%%%%%%
%%%%%%%%%%%%
\section{Classical analytic-type parabolic estimates}
\label{estimthetahorizon}\setcounter{equation}{0} 

The purpose of this section is to prove  Proposition\refer{conventionalestimate}.  We shall  use  the algebraic structure of the Navier-Stokes system and the fact that the function~$\Phi$ is subadditive. 

Let us first bound the horizontal component. We recall that
$$
w_\Phi^h(t) =e^{t\D_\e +\Phi(t,D)} w^h(0) -\bigl(\EE M^\perp(v\cdot \nabla w^h)\bigr)_\Phi (t)
-\bigl(\EE M^\perp ( w^3\partial_3 \vbar^h)\bigr)_\Phi(t) -\bigl( \EE(\nabla_h q)\bigr)_\Phi(t).
$$
We note that $v\cdot \nabla w^h=\dive_h(v^h \otimes w^h)+\partial_3(w^3 w^h)$, recalling
that~$v^{3} = w^3$. On the one hand, using Lemma\refer{PhiactonEE} and Corollary\refer{paradiffanalytic}, we can 
write 
\begin{eqnarray*}
\e\|\EE(\dive_h(v^hw^h))_\Phi\|_{L^1_T(B^{\frac 72})}&\leq&
 C\e\|(v^h w^h)_\Phi\|_{L^1_T(B^{\frac 72})} \\
&\leq &C\|v_\Phi\|_{\wt L^\infty_T(B^{\frac  72})}\e\|w^h\|_{L^1_T(B^{\frac 72})}.
\end{eqnarray*}
By definition of $\theta$,  we infer that
\beq
\label{democonventionalestimateeq00}
\e\|\EE(\dive_h(v^hw^h))_\Phi\|_{L^1_T(B^{\frac 72})}\leq C\theta(T)\|v^h_\Phi\|_{\wt L^\infty_T(B^{\frac 72})}.
\eeq
 On the other hand,  Lemma\refer{PhiactonEEepsd3} and Corollary\refer{paradiffanalytic} imply that
\ben
 \|\EE\big( \e\partial_3 M^\perp(w^3w^h)\big)_\Phi\|_{L^1_T(B^{\frac 72})} & \leq &  C\|w^3w^h\|_{L^1_T(B^{\frac 72})}\nonumber\\
 \label{democonventionalestimateeq01}
 & \leq & C\theta(T)\|v^h_\Phi\|_{\wt L^\infty_T(B^{\frac 72})}.
\een
 For the second term, we use paradifferential calculus which gives 
 \begin{eqnarray*}
 w^3\partial_3\vbar^h&=&T_{w^3}\partial_3\vbar^h+R_{\partial_3\vbar^h} w^3\\
&=&\partial_3T_{w^3}\vbar^h - T_{\partial_3 w^3}\vbar^h+R_{\partial_3\vbar^h} w^3.
 \end{eqnarray*}
  Using again Lemma\refer{PhiactonEEepsd3} and Corollary\refer{paradiffanalytic},  we get
\beno
\|\EE\big( \e\partial_3 M^\perp T_{w^3}\vbar^h)_\Phi\|_{L^1_T(B^{\frac 72})}
& \leq & C\|(T_{w^3}\vbar^h)_\Phi\|_{L^1_T(B^{\frac 72})}\\
 & \leq & C\|w^3_\Phi\|_{ L^1_T(B^{\frac  72})}\|\vbar^h\|_{\wt L^\infty_T(B^{\frac 72})}.
\eeno
By definition of~$\theta$, we infer
\beq
\label{democonventionalestimateeq02}
\|\EE\big( \e\partial_3 M^\perp T_{w^3}\vbar^h)_\Phi\|_{L^1_T(B^{\frac 72})}\leq C\theta(T)\|v^h_\Phi\|_{\wt L^\infty_T(B^{\frac 72})}.
\eeq
By Lemma\refer{PhiactonEE} and Corollary\refer{paradiffanalytic}, we can write that
\beno
\|\EE\big( \e M^\perp T_{\partial_{3}w^3}\vbar^h)_\Phi\|_{L^1_T(B^{\frac 72})}
& \leq & C \e \|(T_{\partial_{3}w^3}\vbar^h)_\Phi\|_{L^1_T(B^{\frac 72})}\\
 & \leq & C \e \|w^3_\Phi\|_{ L^1_T(B^{\frac  72})}
\|\vbar^h\|_{\wt L^\infty_T(B^{\frac 72})}
\eeno
so that
\beq
\label{democonventionalestimateeq03}
\|\EE\big( \e M^\perp T_{\partial_{3}w^3}\vbar^h)_\Phi\|_{L^1_T(B^{\frac 72})}\leq
 C\theta(T)\|v^h_\Phi\|_{\wt L^\infty_T(B^{\frac 72})},
\eeq
and finally along the same lines we have
\beq
\label{democonventionalestimateeq04}
\|\EE\big( \e M^\perp R_{\partial_{3}\vbar^h}w^{3})_\Phi\|_{L^1_T(B^{\frac 72})}\leq
 C\theta(T)\|v^h_\Phi\|_{\wt L^\infty_T(B^{\frac 72})}.
\eeq
 Now we are left with the study of the pressure. Some of its properties are described  in the following lemma.
\begin{lemma}
\label{estimpressureL1}
{\sl  Let us define~$\nabla_\e\eqdefa (\nabla_h,\e\partial_3)$. The following
 two inequalities on the rescaled pressure hold:
\beno
\e\|(\EE \nabla_\e M^\perp q )_\Phi\|_{\wt L^\infty_T(B^{\frac 72})}  & \leq &  C\|v_\Phi\|^2_{\wt L^\infty_T(B^{\frac 72})}\andf \\
\e\|(\EE \nabla_\e M^\perp q )_\Phi\|_{L^1_T(B^{\frac 72})}  & \leq & C\|v_\Phi\|_{\wt L^\infty_T(B^{\frac 72})}\theta(T) .
\eeno
}
\end{lemma}
\begin{proof}
Using the formula\refeq{eqrescaledpressure} on the rescaled pressure and the divergence free condition on~$v$, let us decompose it as~$\e q=q_{1,\e}-q_{2,\e}$ with
\beno
q_{1,\e} & \eqdefa  & \sum_{ k=1}^2 \partial_k\partial_\ell \Delta_\e^{-1}(\e w^kv^\ell)
+ \sum\limits_{1\leq k\leq 2}\partial_k(\e\partial_3)\Delta_\e^{-1}(w^3 v^k)\andf\\
q_{2,\e} & \eqdefa &2\e\partial_3\Delta_\e^{-1}(w^3\dive_h w^h).
\eeno
Let us start with~$q_{1,\e}$. We have
$$
\nabla_\e q_{1,\e} = \sum_{k=1}^2 \partial_k\biggl(  \sum_{\ell=1}^2 \nabla_\e\partial_\ell\Delta_\e^{-1}(\e w^k v^\ell) + \nabla_\e (\e\partial_3) \D_\e^{-1}(w^3v^k)\biggr).
$$
As~$\nabla_\e^2\D_\e^{-1}$ is a family of bounded Fourier multipliers (uniformly with respect
 to~$\e$), we infer from Lemma\refer{PhiactonEE} and Corollary\refer{paradiffanalytic} 
that
\ben
\label{estimpressureL1demoeq1}
\e\|(\EE(\nabla_\e M^\perp q_{1,\e}))_\Phi\|_{\wt L^\infty_T(B^{\frac 72})}  & \leq &  C\|v_\Phi\|^2_{\wt L^\infty_T(B^{\frac 72})}\andf \\
\label{estimpressureL1demoeq2}
\e\|(\EE(\nabla_\e M^\perp q_{1,\e}))_\Phi\|_{L^1_T(B^{\frac 72})}  & \leq & C\|v_\Phi\|_{\wt L^\infty_T(B^{\frac 72})}\theta(T) .
\een
In order to study~$q_{2,\e}$, let us observe that
\ben
w^3\dive_h w^h & = & R_{\dive_h w^h} w^3+ T_{w^3} \dive_h w^h\nonumber \\
\label{estimpressureL1demoeq3}
 & =& R_{\dive_h w^h} w^3 +\sum_{k=1} ^2 \bigl(\partial_kT_{w^3} w^k - T_{\partial_kw^3} w^k\bigr).
\een
As above we get,  using Lemma\refer{PhiactonEE} and 
Corollary\refer{paradiffanalytic},
\beno
\e\|(\EE(\nabla_\e M^\perp q_{2,\e}))_\Phi\|_{\wt L^\infty_T(B^{\frac 72})}  & \leq &  C\|v_\Phi\|^2_{\wt L^\infty_T(B^{\frac 72})}\andf \\
\e\|(\EE(\nabla_\e M^\perp q_{2,\e}))_\Phi\|_{L^1_T(B^{\frac 72})}  & \leq & C\|v_\Phi\|_{\wt L^\infty_T(B^{\frac 72})}\theta(T) .
\eeno
Together with estimates\refeq{estimpressureL1demoeq1} and\refeq{estimpressureL1demoeq2}, this concludes the proof of the lemma.
\end{proof}

The above Lemma\refer{estimpressureL1},
 together with estimates\refeq{democonventionalestimateeq00} to\refeq{democonventionalestimateeq03},
 implies that
\beq
\label{conventionaldemoeq1}
\e \|w^h\|_{L^1_T(B^{\frac 7 2 })} \leq \e \|e^{a|D_3|}w^h_0\|_{B^{\frac 7 2}} + C_0^{(1)} \|v_\Phi\|_{\wt L^\infty_T(B^{\frac 7 2 })}\displaystyle  \theta (T) .
\eeq

Let us prove the estimates on the vertical component. It turns out that  it is better behaved because of the special structure of the system. Indeed, 
thanks to the divergence free  condition, almost no vertical derivatives appear in the equation of~$w^3$: we have (since~$w^{3} = v^{3}$)
\beq 
\label{eqverticalcomponent}
\partial_t w^3 -\D_\e  w^3 = -v^h\cdot\nabla_h w^3+w^3\dive _h w^h -\e^2\partial_3 q.
\eeq
The Duhamel formula reads
$$
w^3(t) =e^{t\D_\e} w^3(0) + \EE M^\perp (w^3\dive_h w^h-v^h\cdot \nabla_h w^3)(t) 
-\EE M^\perp (\e^2\partial_3q)(t) .
$$
Applying the Fourier multiplier~$e^{\Phi(t,D)}$ to the above relation gives
\beq
\label{DuhamePhivertical}
w_\Phi^3(t) =e^{t\D_\e +\Phi(t,D)} w^3(0) + \bigl(\EE M^\perp(w^3\dive_h w^h-v^h\cdot \nabla_hw^3)\bigr)_\Phi(t) 
-\bigl(\EE M^\perp \e^2\partial_3q\bigr)_\Phi(t) .
\eeq
Using\refeq{estimpressureL1demoeq3} and then Lemma\refer{PhiactonEE} and Corollary\refer{paradiffanalytic}, we get
\ben
\label{conventionalestimatedemoeq1}\bigl\|\bigl(\EE M^\perp(w^3\dive_h w^h)\bigr)_\Phi\bigr\|_{\wt L^{\infty}_{T}(B^{\frac72})} 
& \leq & C
\|w^{3}_{\Phi}\|_{\wt L^{\infty}_{T}(B^{\frac72})} \|w^{h}_{\Phi}\|_{\wt L^{\infty}_{T}(B^{\frac72})}  \andf\\
\label{conventionalestimatedemoeq2} \bigl\|\bigl(\EE M^\perp(w^3\dive_h w^h)\bigr)_\Phi\|_{ L^{1}_{T}(B^{\frac72})} 
& \leq & C
\|w^{3}_{\Phi}\|_{L^{1}_{T}(B^{\frac72})} \|w^{h}_{\Phi}\|_{\wt L^{\infty}_{T}(B^{\frac72})}.
\een
Writing that
\beno
v^h\cdot\nabla_h a & = & \sum_{k=1}^2 \bigl(T_{v^k}\partial_k a+R_{\partial_k a} v^k\bigr)\nonumber\\
 &= & -T _{\dive_h w^h}  a + \sum_{k=1}^2 \bigl(\partial_k T_{v^k} a+R_{\partial_k a} v^k\bigr)
\eeno
and using Lemma\refer{PhiactonEE} and Corollary\refer{paradiffanalytic}, we get
\beno
\bigl\|\bigl(\EE M^\perp(v^h\cdot\nabla_h w^3)\bigr)_\Phi\bigr\|_{\wt L^{\infty}_{T}(B^{{\frac72}})} 
& \leq & C
\|w^{3}_{\Phi}\|_{\wt L^{\infty}_{T}(B^{{\frac72}})} \|w^{h}_{\Phi}\|_{\wt L^{\infty}_{T}(B^{{\frac72}})}  \andf\\
\bigl\|\bigl(\EE M^\perp(v^h\cdot\nabla_h w^3)\bigr)_\Phi\|_{ L^{1}_{T}(B^{{\frac72}})} 
& \leq & C
\|w^{3}_{\Phi}\|_{L^{1}_{T}(B^{{\frac72}})} \|w^{h}_{\Phi}\|_{\wt L^{\infty}_{T}(B^{{\frac72}})}.
\eeno
Together with estimates\refeq{conventionalestimatedemoeq1} and\refeq{conventionalestimatedemoeq2}, and Lemma\refer{estimpressureL1}, this gives
\beno
\|w^3\|_{L^1_T(B^{\frac 7 2})}  & \leq &  \|e^{a|D_3|}w^3_0\|_{B^{\frac 7 2}} + C_0^{(1)} \|v_\Phi\|_{\wt L^\infty_T(B^{\frac 7 2})}  \theta (T)\andf \\
\|w^3_\Phi\|_{\wt L^\infty_T(B^{\frac 7 2})} &  \leq &  \|e^{a|D_3|}w^3_0\|_{B^{\frac 7 2}} + C_0^{(1)} \|v_\Phi\|_{\wt L^\infty_T(B^{\frac 7 2})}^2.
\eeno
Together with\refeq {conventionaldemoeq1}, this concludes the proof of Proposition\refer{conventionalestimate}.

%%%%%%%%%%%%
%%%%%%%%%%%%

\section{The gain of one vertical derivative on the horizontal part}
\setcounter{equation}{0} 
In this section we shall prove Proposition~\ref{estimhorinzonLinftytilde}. The proof will be separated into two parts: first we shall consider the case of the horizontal average~$\vbar_\Phi^h$, and then the remainder~$w_\Phi^h$.

\subsection{The gain of one vertical derivative on the horizontal average}
\label{thetavbargain}
We shall study in this section the equation on the horizontal average of the solution. We emphasize that in the equation on $\vbar$ we cannot recover the vertical derivative appearing in the force term by the regularizing effect.  The fundamental idea to gain a vertical derivative is to  use the analyticity of the solution and therefore to estimate $\vbar_\Phi$. The lemma is the following.

\begin{lemma}
\label{lemmeestivbar}
{\sl A constant~$C_0 $  exists such that, for any positive~$\lam$,  for any initial data~$v_0$, and 
 for any~$T$ satisfying~$\displaystyle   \theta (T) \leq  a / \lam $, we have
$$
\|\vbar_\Phi^h\|_{\wt L^\infty_T(B^\frac72)} \leq  \|e^{a|D_3|} \vbar^h_0\|_{B^\frac72}
 + C_0 \Bigl( \frac {1} \lam  +  \|v_\Phi\|_{\wt L^\infty_T(B^\frac72)}\Bigr)
 \|v^h_\Phi\|_{ \wt L^\infty_T(B^\frac72)}.
$$
}
\end{lemma}

\begin{proof}
The horizontal average~$\vbar$ satisfies
\beq
\label {lemmeestivbardemoeq0}
\partial_{t} \vbar -\e^2\partial_3^2 \vbar = -\partial_3 M (w^3w^h)\andf
\vbar_{|t=0}=\vbar_{0}.
\eeq
Let us define~$G\eqdefa  -\partial_3 M ( w^3w^k)$. Writing the solution of\refeq {lemmeestivbardemoeq0} in terms of the Fourier transform, we get, using\refeq {relationLerayNantes} with~$\xi_h=0$,
$$
|\cF (\vbar_\Phi) (t,\xi) | \leq  |\cF\vbar_0(\xi)| e^{a|\xi_3|}
+\int_0^t e^{-\lam |\xi_3|\int_{t'} ^t \dot \theta (t'')dt''} |\cF( G_\Phi) (t',\xi)|dt'.
$$
Then,  taking the~$L_j^2$ norm, we infer that
\beq
\label{lemmeestivbardemoeq01}
\|\vbar_\Phi(t) \|_{L_j^2} \leq \|e^{a|D_3|} \vbar_0\|_{L_j^2} 
+ \int_0^t  e^{- c\lam2^j \int_{t'} ^t \dot \theta (t'')dt''} \| G_\Phi(t')\|_{L_j^2} dt'.
\eeq
Now, let us estimate~$ \| G_\Phi(t')\|_{L_j^2}$.
For any function~$a$, using the fact that the vector field~$w$ is divergence free, let us write that
\ben
\partial_3 (w^3 a) & = & \partial_3 \big ( T_{w^3} a +R_{a} w^3\bigr)\nonumber \\
 & = &  \partial_3 T_{w^3} a  +R_{\partial_3 a } w^3 -R _{a} \dive_h w^h\nonumber\\
\label{lemmeestivbardemoeq000}& = & \partial_3 T_{w^3} a  
+R_{\partial_3 a} w^3 -\sum_{\ell=1}^2 \partial_\ell  R _{a} w^\ell  +\sum_{\ell=1}^2 R_{\partial_\ell a} w^\ell.
\een
Thus, we infer that
\ben
G &  =  & - \partial_3 MT_{w^3} w^k 
-M \Bigl( R_{\partial_3 w^k} w^3 +\sum_{\ell=1}^2 R_{\partial_\ell w^k} w^\ell
 -\sum_{\ell=1}^2 \partial_\ell  R _{w^{3}} w^\ell \Bigr)\nonumber\\
 \label{lemmeestivbardemoeq1}
 & = &   -\partial_3 MT_{w^3} w^k 
-M \Bigl( R_{\partial_3 w^k} w^3 +\sum_{\ell=1}^2 R_{\partial_\ell w^k} w^\ell\Bigr).
\een
Now, let us study $ \cF M (T_a b)_\Phi $ and~$\cF M (R_a b)_\Phi$ 
for two functions $a$ and $b$ which have~$0$ horizontal average. As the two terms are identical, let us study the first one. By definition, we have
$$
\cF \bigl ( T_a b ) (t,(0,\xi_3))  =   \sum_{j}\int_{2^{j}\cC \cap B((0,\xi_3), 2^{j}) } \!\!\wh a((0,\xi_3)-\eta) \wh b(\eta) d\eta.
$$
As~$\theta(T)\leq \lam^{-1}a$,  by definition of~$\Phi$  we have, for any~Ê$\eta \in (\ZZ^2\setminus\{0\})\times\R$,
\beno
\Phi (t,(0,\xi_3)) & \leq & \Phi  (t,(0,\xi_3-\eta_3)) + \Phi  (t,(0,\eta_3))\\
 & \leq & -2t^{\frac 1 2} +\Phi( t, ((0,\xi_3)-\eta) +\Phi (t, -\eta).
 \eeno
Thus we have
$$
| (\cF M (T_a b)_\Phi ) (t,\xi) |  \leq   
e^{-2t^{\frac 1 2} } (\cF M T_{a^+_\Phi} b^+_\Phi ) (t,\xi).
$$
Applied to\refeq{lemmeestivbardemoeq1}, this implies that
$$\longformule{
\bigl| \cF G_\Phi (t,\xi)\bigr| \leq
|\xi_3| \cF \bigl(T_{(w^3_\Phi)^+} (w^k_\Phi)^+\bigr) (t,(0,\xi_3))}
{
{}+ e^{-2t^{\frac 12}} 
\cF \Bigl( R_{(\partial_3 w^k_\Phi) ^+} (w^3_\Phi ) ^+ +\sum_{\ell=1}^2 
R_{(\partial_\ell w^k_\Phi)^+} (w^\ell_\Phi)^+\Bigr) (t,(0,\xi_3)).
}$$
Inequality\refeq{estmitildepastilde} then implies that, for any~$t\in [0,T]$, 
\beq
\label{lemmeestivbardemoeq2}
2^{j\frac 7 2 } \| G_\Phi (t)\|_{L_j^2} \leq
C c_j  \|v^h_\Phi\|_{\wt L^\infty_T(B^{\frac 7 2})} \bigl( 2^j \|w_\Phi^3(t)\|_{B^{\frac 7 2}} 
+e^{-2t^{\frac 1 2} }  \|v_\Phi\|_{\wt L^\infty_T(B^{\frac 7 2})}\bigr).
\eeq

Then, by definition of~$\theta$, Inequalities\refeq{lemmeestivbardemoeq01} and\refeq{lemmeestivbardemoeq2} imply that 
$$
\longformule{
2^{j{\frac 7 2}} \|(\vbar_\Phi)(t) \|_{L_j^2} \leq 2^{j{\frac 7 2}} \|e^{a|D_3|} \vbar_0\|_{L_j^2} 
}
{
{}+
C c_j \|v_\Phi\|_{\wt L^\infty_T(B^{\frac 7 2})} \biggl(\int_0^t e^{-c\lam 2^j 
\int_{t'} ^t \dot \theta (t'')dt''}
2^j\dot \theta(t')dt'+ \|v_\Phi\|_{\wt L^\infty_T(B^{\frac 7 2})}  \int_0^t e^{-2t'^{\frac 1 2} }  dt'\biggr).
}
$$
This gives
$$
2^{j{\frac 7 2}} \|\vbar_\Phi \|_{L^\infty_T(L_j^2)} \leq 2^{j{\frac 7 2}} \|e^{a|D_3|} \vbar_0\|_{L_j^2} 
+ C c_j \|v^h_\Phi\|_{\wt L^\infty_T(B^{\frac 7 2})} \Bigl(\frac 1 \lam + \|v_\Phi\|_{\wt L^\infty_T(B^{\frac 7 2})} \Bigr).
$$
Taking the sum over~$j$ concludes the proof of the lemma.
\end{proof}

%%%%%%%%%%%%%%%%%%%%%%%%%%%%%%%%%%%%%%%%%%%%%%%%%

\subsection{The gain of the vertical derivative on the whole horizontal term}
\label{thegain}
Now let us estimate the rest of the horizontal  term, that is~$ \|w^h_\Phi\|_{\wt L^\infty_T(B^{\frac 7 2})}$. 
As in Section\refer{thetavbargain}, the function~$\theta$ will play a crucial role.
\begin{lemma}
\label{estimLinftyhorinzon}
{\sl A constant~$C_0 $ exists such that, for any~$\lam$,  for any initial data~$v_0$, and 
 for any~$T$ satisfying~$\displaystyle   \theta (T) \leq  a / \lam $, we have
$$
\|w^h_\Phi\|_{\wt L^\infty_T(B^{\frac 7 2})} \leq \|e^{a|D_3|}w^h _0\|_{B^{\frac 7 2}}  +
C_0  \Bigl( \frac {1} \lam  +  \|v_\Phi\|_{\wt L^\infty_T(B^{\frac 7 2})}\Bigr) \|v^h_\Phi\|_{ \wt L^\infty_T(B^{\frac 7 2})}.
$$
}
\end{lemma}
\begin{proof}
The Duhamel formula writes
$$ 
w^h(t) =e^{t\D_\e} w^h(0)  -\EE\dive_h (v^h\otimes v^h   )(t)
 -\EE M^\perp\partial_3 ( w^3 v^h )(t)- \EE(\nabla_h q)(t). 
$$
Lemma\refer{PhiactonEE} and Corollary\refer {paradiffanalytic} imply that
\ben
\|(\EE\dive_h (v^h\otimes v^h ))_\Phi\|_{\wt L^\infty_T(B^{\frac 7 2})} &  \leq  &
C  \|  (v^h\otimes v^h  )_\Phi\|_{\wt L^\infty_T(B^{\frac 7 2})} \nonumber \\
\label{estimLinftyhorinzondemoeq1}
& \leq  &  C \| v^h_\Phi\|^2_{\wt L^\infty_T(B^{\frac 7 2})}.
\een
Then  using\refeq{lemmeestivbardemoeq000} we get,  thanks  to Leibnitz formula,
\beno
M^\perp \partial_3 (w^3 v^k) & = & M^\perp \partial_3 T_{w_3} v^k +F^k \with\\
 F^k & \eqdefa &    M^\perp \biggl(R_{\partial_3v^k} w^3-\sum_{\ell=1}^2  \bigl(\partial_\ell R_{v^k} w^\ell -R_{\partial_\ell v_k} w^\ell\bigr) \biggr).
\eeno
Thanks to Lemma\refer{PhiactonEE} and Corollary\refer {paradiffanalytic}, we get
\beno
\| (\EE M^\perp F^k)_\Phi\| _{\wt L^\infty_T(B^{\frac 7 2})}  & \leq & 
\|(F^k)_\Phi\| _{\wt L^\infty_T(B^{\frac 7 2})}\\
& \leq & C\|v_\Phi\| _{\wt L^\infty_T(B^{\frac 7 2})}\|v^h_\Phi\| _{\wt L^\infty_T(B^{\frac 7 2})}.
\eeno
Together with Lemma\refer{PhiactonEEd3gainlambda}, this gives
\beq
\label{estimLinftyhorinzondemoeq111}
\bigl\|\bigl(\EE M^\perp\partial_3 ( w^3 v^h )\bigr)_\Phi \bigr\|_{\wt L^\infty_T(B^{\frac 7 2})} \leq 
C_0  \Bigl( \frac {1} \lam  +  \|v_\Phi\|_{\wt L^\infty_T(B^{\frac 7 2})}\Bigr) \|v^h_\Phi\|_{ \wt L^\infty_T(B^{\frac 7 2})}.
\eeq
Now let us study the pressure term.  Formula\refeq{eqrescaledpressure} together with
 the divergence free condition leads to the decomposition~$q=q_h(w^h)+q_3(v)$ with
\begin{eqnarray}
\label{definpressh}
\quad\quad q_h(w^h) & \eqdefa &  \D_\e^{-1}\Bigl((\dive_h w^h)^2+ \sum_{1\leq k,\ell\leq 2} \partial_k w^\ell
\partial_\ell w^k\Bigr)\andf \\
\label{definpress3}
q_3(v) & \eqdefa  & \D_\e^{-1}\Bigl( \sum_{1\leq \ell\leq 2}  \partial_3 v^\ell
\partial_\ell w^3\Bigr).
\end{eqnarray}
For the first term we use    Bony's decomposition in order to obtain
$$
\partial_k w^\ell\partial_\ell w^k = T_{\partial_k w^\ell} \partial_\ell w^k +R_{\partial_\ell w^k}\partial_k w^\ell.
$$
Then the Leibnitz formula  implies that
\begin{equation}
\label{paradiff:h1}
\partial_k w^\ell\partial_\ell w^k =\partial_\ell T_{\partial_k w^\ell}  w^k + \partial_k R_{\partial_\ell w^k}w^\ell -  T_{\partial_\ell \partial_k w^\ell}  w^k - R_{ \partial_k\partial_\ell w^k}w^\ell.
\end{equation}
On the other hand, again by   paradifferential calculus, we can write that
\begin{eqnarray}
\label{paradiff:h2}
\nonumber(\dive_h w^h)^2  & = & T_{\dive_h w^h} \dive_h w^h +R_{\dive_h w^h} \dive_h w^h\\
 & = & \dive_h \bigl( T_{\dive_h w^h} w^h + R_{\dive_h w^h} w^h \bigr)
 -\sum_{k=1}^2 \left(T_{\partial_k \dive_h w^h} w^k +R_{\partial_k \dive_h w^h} w^k\right).
\end{eqnarray}
Then Lemma\refer{PhiactonEE} implies that
$$
\|(\EE\nabla_h q_h(w^h))_\Phi\|_{\wt L^\infty_T(B^{\frac 72})} \leq  C_0 \|(M^\perp q_h(w^h))_\Phi\|_{\wt L^\infty(B^{\frac 72})}.
$$
Using Corollary \ref{paradiffanalytic} and the fact that the operators
 $\nabla_h\Delta_\e^{-1}M^\perp$ and $\Delta_\e^{-1} M^{\perp}$ are bounded (uniformly
in~$\e$) Fourier multipliers,  we obtain  
\beq\label{estimatefinalqhwh}
\bigl\|\bigl( \EE(\nabla_h q_h(w^h))\bigr)_\Phi\bigr\|_{\wt L^\infty_T(B^{\frac 72})} \leq C_0 \|v^h_\Phi\|_{\wt L^\infty_T(B^{\frac 72})}^2.
\eeq
For the second term, let us decompose~$q_3(v)$ in the following way:
\begin{eqnarray*}
\partial_3 v^\ell\partial_\ell w^3&=& T_{\partial_3 v^\ell} \partial_\ell w^3
+ R_{\partial_\ell w^3}\partial_3 v^\ell\\
&=&\partial_\ell T_{\partial_3 v^\ell} w^3+\partial_3 R_{\partial_\ell w^3} v^\ell-T_{\partial_3\partial_\ell v^\ell} w^3-R_{\partial_3\partial_\ell w^3} v^\ell.
\end{eqnarray*}
Using now Lemma \ref{PhiactonEE} together with Corollary\refer{paradiffanalytic} and Lemma\refer{PhiactonEEd3gainlambda}, we obtain
\beq\label{estimatefinalq3w3}
\bigl\|\bigl( \EE(\nabla_h q_3(v))\bigr)_\Phi\bigr\|_{\wt L^\infty_T(B^{\frac 72})}\leq 
C_0 \Bigl(\frac{1}{\lam}+\|v_\Phi\|_{\wt L^\infty_T(B^{\frac 72})}\Bigr)\|v^h\|_{\wt L^\infty_T(B^{\frac 72})}.
\eeq
The expected result is obtained putting together estimates~(\ref{estimatefinalqhwh}) and\refeq{estimatefinalq3w3}
on the pressure with estimates\refeq{estimLinftyhorinzondemoeq1}  and\refeq{estimLinftyhorinzondemoeq111}  on the nonlinear terms.
\end{proof}

%%%%%%%%%%%%%%%%%%%%%%%%%%%%%%%%%%%%%%%%%%%%%%%%%%%%%%%%%%%%%%%
%%% BIBLIO
%%%%%%%%%%%%%%%%%%%%%%%%%%%%%%%%%%%%%%%%%%%%%%%%%%%%%%%%%%%%%
%{\tt v\'erifier biblio}


\begin{thebibliography}{99} 

\bibitem{cannonemeyerplanchon}
M. Cannone, Y. Meyer and F. Planchon,
Solutions autosimilaires des {\'e}quations de Navier-Stokes,
{S{\'e}minaire "{\'E}quations aux D{\'e}riv{\'e}es Partielles" de l'{\'E}cole polytechnique},
Expos{\'e} VIII, 1993--1994.

\bibitem{chemin21}
J.-Y. Chemin,
Le syst\`eme de Navier-Stokes incompressible soixante dix
ans apr\`es Jean Leray,
S\'eminaire et Congr\`es, {\bf 9}, 2004,   pages   99--123.

\bibitem{cdgg2}
J.-Y. Chemin, B. Desjardins, I. Gallagher and E. Grenier,
Fluids with anisotropic viscosity,
{\em Mod\'elisation Math\'ematique et Analyse Num\'erique}, {\bf 34}, 2000,
pages~315--335.

\bibitem{cgens} 
J.-Y. Chemin and I. Gallagher, On  the global wellposedness  of the 3-D 
   Navier-Stokes equations with large initial data, {\it Annales de l'\'Ecole Normale  Sup\'erieure}, {\bf 39}, 2006, pages 679--698.    

\bibitem{cg2} 
J.-Y. Chemin and I. Gallagher, Wellposedness and stability results for  the  Navier-Stokes equations 
in~${\bf R}^{3}$
{\it to appear in    Annales de l'Institut Henri Poincar\'e, Analyse Non Lin\'eaire.
}

\bibitem{cg3}
J.-Y. Chemin and I. Gallagher,
Large, global solutions to the Navier-Stokes  equations,   slowly varying in one direction,
 {\it to appear in Transactions of the Americal Mathematical Society}.
 
 \bibitem{chemin13}
J.-Y. Chemin and N. Lerner,
Flot de champs de vecteurs non-lipschitziens et \'equations de Navier-Stokes,
{\em Journal of Differential Equations}, {\bf 121}, 1995, pages 314--328.

\bibitem{cz1} 
J.-Y. Chemin and P. Zhang, On the global  wellposedness of the 3-D incompressible anisotropic
Navier-Stokes equations,
{\it Communications in  Mathematical Physics}, {\bf 272}, 2007, pages 529--566.

\bibitem{ESS}
L. Escauriaza, G. Seregin and V. \c{S}ver\'ak, On the $L_{3,\infty}$ solutions to the Navier-
Stokes equations and backward uniqueness, Russ. Math. Surv., {\bf 58}, 2003,  pages  211--250.

\bibitem{fujitakato}
H. Fujita and T. Kato,
On the Navier-Stokes initial value problem I,
{\em  Archive for Rational Mechanics and Analysis}, 
{\bf 16}, 1964, pages~269--315.

\bibitem{gallagheriftimieplanchon} I. Gallagher, D. Iftimie and F. Planchon, Asymptotics
   and stability for global solutions to the
   Navier--Stokes equations, {\it  Annales de l'Institut Fourier},  {\bf 53}, 2003, pages 1387--1424. 

\bibitem{gp}  I. Gallagher and M. Paicu, Remarks on the blow-up  of solutions to  a   toy model for the Navier-Stokes equations, {\it submitted}.

\bibitem{gigamiyakawa} Y.   Giga  and T. Miyakawa, 
Solutions in $L^r$ of the Navier-Stokes initial value problem, {\em Archiv for Rational Mechanics and  Analysis}, {\bf 89}, 1985, pages 267--281.

\bibitem{dragosunique}
D. Iftimie, A uniqueness result for the Navier-Stokes equations with vanishing vertical viscosity, {\it SIAM Journal of  Mathematical  Analysis}, {\bf  33}, 2002), pages 1483--1493.

\bibitem{iftimieraugelsell} D. Iftimie, G. Raugel and G.R. Sell, 
Navier-Stokes equations in thin 3D domains with Navier boundary conditions, {\em Indiana University Mathematical Journal}, {\bf 56}, 2007, pages 1083--1156.

\bibitem{kato} T. Kato,
Strong $L^p$-solutions of the Navier-Stokes equation in
$\R^m$ with applications to weak solutions, {\it Mathematische
Zeitschrift}, {\bf 187}, 1984,  pages 471-480 .

\bibitem{kochtataru} H. Koch and D. Tataru, Well--posedness for the
   Navier--Stokes equations, {\it Advances in
Mathematics},  {\bf 157}, 2001, pages~22--35.

\bibitem{ladyzhenskaya} O. Ladyzhenskaya,  {\it 
 The mathematical theory of viscous incompressible flow}. Second English edition, revised and enlarged. T Mathematics and its Applications, Vol. 2 Gordon and Breach, Science Publishers, New York-London-Paris 1969 xviii+224 pp.


\bibitem{mahalovtitileibovich} S.  Leibovich, A. Mahalov, E.  Titi,  
Invariant helical subspaces for the Navier-Stokes equations. 
{\em  Archiv for Rational Mechanics and  Analysis},  {\bf 112}, 1990,  pages 193--222. 



\bibitem{lemarie} P.-G. Lemari\'e-Rieusset,  Recent developments in the Navier-Stokes problem. Chapman and Hall/CRC Research Notes in Mathematics, {\bf 43}, 2002.

\bibitem{leray}
J. Leray, Essai sur le mouvement d'un liquide visqueux emplissant l'espace,
{\em Acta Matematica}, {\bf 63}, 1933, pages 193--248.

\bibitem{leray2D} J.~Leray,
\'Etude de diverses \'equations int\'egrales non lin\'eaires et
de quelques probl\`emes que pose l'hydrodynamique.
{\em Journal de  Math\'ematiques Pures et  Appliqu\'ees}, {\bf 12}, 1933,  pages 1--82. 

\bibitem{babinandco} A. Mahalov and B. Nicolaenko, 
 Global solvability of three-dimensional Navier-Stokes equations with
 uniformly high initial vorticity, 
  (Russian. Russian summary)
{\em Uspekhi Mat. Nauk} {\bf 58}, 2003, pages 79--110; translation in
 {\em Russian Math. Surveys} {\bf 58}, 2003, pages 287--318.


\bibitem{meyer} Y. Meyer, 
{\sl Wavelets, Paraproducts and Navier--Stokes}.  
Current Developments in Mathematics, International Press, Cambridge, 
Massachussets, 1996.

\bibitem{ms}
S. Montgomery-Smith, Finite-time blow up for a Navier-Stokes like equation,  {\it Proc. Amer. Math. Soc.} {\bf 129}, 2001,  pages 3025--3029.


\bibitem{NRS}
J. Ne\c{c}as, M. Ru\c{z}i\c{c}ka and V. \c{S}ver\'ak, On Leray's self-similar solutions of the
Navier-Stokes equations, {\it Acta Mathematica}, {\bf 176}, 1996,  pages  283--294.

\bibitem{marius1}
M.Paicu, \'Equation anisotrope de Navier-Stokes dans des espaces critiques,  {\it Revista Rev. Matemtica Iberoamericana}, {\bf  21},  2005, page 179--235. 

\bibitem{poncerackesideristiti} G. Ponce, R.  Racke, T.  Sideris and E.  Titi,  Global stability of large solutions to the $3$D Navier-Stokes equations, {\em Communications in  Mathematical  Physics},  {\bf 159}, 1994,  pages 329--341.

\bibitem{raugelsell}
G. Raugel and G.R. Sell,  Navier-Stokes equations on thin $3$D domains. I.
 Global attractors and global regularity of solutions,
{\it   Journal of the American  Mathematical  Society},~{\bf 6}, 1993,
pages 503--568.

\bibitem{sammartino&caflisch}
M. Sammartino and R. E. Caflisch,  Zero Viscosity Limit for Analytic Solutions, of the
Navier-Stokes Equation on a Half-Space. I. Existence for Euler and Prandtl Equations, {\em Communications in Mathematical Physics}, {\bf 192}, 1998,  pages 433--461.

 \bibitem{ukhoiudo} M. Ukhovskii and V. Iudovich,  Axially symmetric flows of ideal and viscous fluids filling the whole space. {\em Prikl. Mat. Meh.} 32 59--69 (Russian); translated as {\em 
 Journal of  Applied  Mathematics and  Mechanics},  {\bf 32}, 1968, pages  52--61. 

\bibitem{weisslerNS}
F. Weissler, The Navier-Stokes Initial Value Problem in~$L^p$, {\it Archiv for Rational Mechanics and Analysis}, {\bf 74}, 1980, pages 219-230.
\end{thebibliography}
\end{document}